\documentclass{amsart}

\usepackage{amssymb}
\usepackage{amsrefs}
\usepackage[all]{xy}

\numberwithin{equation}{subsection}

\newtheorem*{theorem}{Theorem}
\newtheorem*{proposition}{Proposition}
\newtheorem*{QTheorem}{Prequantization Theorem}

\begin{document}

\title{differential characters as stacks and prequantization}

\author{Eugene Lerman}
\author{Anton Malkin}
\thanks{Supported in part by NSF grants DMS-0603892 and DMS-0456714}
\address{Department of Mathematics, 
University of Illinois at Urbana-Champaign, 
1409 W. Green Street, Urbana, IL 61801, USA}

\begin{abstract}
We generalize geometric prequantization of symplectic manifolds 
to differentiable stacks. Our approach is atlas-independent 
and provides a bijection between isomorphism classes of principal
$S^1$-bundles (with or without connections) and second cohomology 
groups of certain chain complexes.
\end{abstract}

\maketitle

\tableofcontents

\section{Introduction.}

\subsection{Quantization and reduction.}

Our original motivation for this project has to do with the 
``quantization commutes with reduction'' principle in 
symplectic geometry. Let us explain the issues involved on
a simple example.

Consider the complex irreducible representations of $SU(2)$.  There is
one in every dimension: $V_0$, $V_1$, \ldots; 
$\dim_{\mathbb{C}} V_n = n$. 
The multiplicity of the zero weight of a maximal torus $T \simeq S^1$ 
of $SU(2)$ in $V_n$ is
\begin{equation}\nonumber
\dim (V_n)^T = \begin{cases}
0 & \text{if $n$ is even;} \\
1 & \text{if $n$ is odd.}
\end{cases} 
\end{equation}
The representation $V_n$ can be constructed by quantizing the complex
projective plane $\mathbb{CP}^1$: 
Start with the pair $(\mathbb{CP}^1, n \omega)$,
where $\omega$ is the $SU(2)$ invariant area form on $\mathbb{CP}^1$
normalized so that $\int_{\mathbb{CP}^1} \omega =1$.  The action of $SU(2)$
on the K\"{a}hler manifold $(\mathbb{CP}^1, n \omega)$ 
is Hamiltonian with an equivariant moment map
$\mu_n: \mathbb{CP}^1\to \mathfrak{su}(2)^*$. 
The K\"{a}hler form $n \omega$ determines a Hermitian holomorphic 
line bundle $L_n \to \mathbb{CP}^1$ with a Hermitian connection.  
The connection and the moment map allow us to lift the action 
of $\mathfrak{su} (2)$ on $\mathbb{CP}^1$ to an action on $L_n$, 
which integrates to an action of
$SU(2)$.  The space of holomorphic sections $H^0(\mathbb{CP}^1, L_n)$ 
is then a representation of $SU(2)$ which happens to be $V_n$, 
and we may think of the Hilbert space $H^0(\mathbb{CP}^1, L_n)$ 
as a quantization of $(\mathbb{CP}^1, n \omega, \mu _n)$.  
By the principle that ``quantization commutes with reduction'' 
\begin{equation}\nonumber
(V_n)^T  = \text{Quant}(\mathbb{CP}^1/ \! / T, (n\omega)_0),
\end{equation}
where $\mathbb{C} P^1/ \! / T$ denotes the symplectic quotient $\mu_n ^{-1}
(0)/T$ and $(n\omega)_0$ the induced symplectic form on the quotient.
It's not hard to see that $\mu_n^{-1} (0)$ is a single $T$-orbit, so
the symplectic quotient $\mathbb{CP}^1 / \! / T$ is a single point.  The
quantization of a point is a complex line bundle over the point, i.e.,
$\mathbb{C} \to \mathrm{pt}$.  
Clearly $H^0 ( \mathrm{pt}, \mathbb{C}) \simeq \mathbb{C}$.  
This contradicts that
$(V_n)^T = 0$ for $n$ even.  What did we do wrong?

Several things.  First of all the symplectic quotient 
$\mathbb{C} P^1/ \! / T$ is not a point.  Rather, it's a point 
with a trivial action of $\mathbb{Z}/2$. We think of it as 
a groupoid $\mathbb{Z}/2 \rightrightarrows \mathrm{pt}$.  A
prequantization of this groupoid is a line bundle over $\mathrm{pt}$ 
with an action of $\mathbb{Z}/2$ and a quantization is the space 
of $\mathbb{Z}/2$-invariant sections of this bundle.  
For $n$ even, the restriction
$L_n|_{\mu_n^{-1} (0)}$ descents to a bundle over $\mathrm{pt}$ 
where $\mathbb{Z}/2$ 
acts non-trivially on the fiber and so the space of invariant sections is
zero.  For $n$ odd, the restriction $L_n |_{\mu_n^{-1} (0)}$ descents to
a bundle over $\mathrm{pt}$ with a trivial $\mathbb{Z}/2$ action and 
so the space of invariant sections is $\mathbb{C}$. 
More abstractly, one can think
of a point with $\mathbb{Z}/2$-action as (an atlas of) a stack 
$[\mathrm{pt} / (\mathbb{Z}/2) ]$. From this point of view the reduction 
procedure is just the restriction of $\omega$ to $\mu_n^{-1} (0)$. Nothing
happens on the quotient stage since $[\mu_n^{-1} (0) / T] =
[\mathrm{pt} / (\mathbb{Z}/2) ]$ as stacks.

Why are there two different ways of (pre)quantizing a point with a
$\mathbb{Z}/2$-action?  After all, there is only one symplectic form on the
groupoid $\mathbb{Z}/2 \rightrightarrows \mathrm{pt}$ (or on 
\mbox{$T \rightrightarrows \mu_n^{-1} (0)$}).  
This only looks puzzling if we think of prequantization as a construction
\begin{equation}\nonumber
\text{integral symplectic form} \mapsto \text{Hermitian line bundle with 
connection}.
\end{equation}
Such a point of view is rather misleading: the
connection is not uniquely determined by its curvature. 
For example flat line bundles are classified by 1-dim representations
of the fundamental group. Nor does
the de Rham cohomology class of the symplectic form see all of the
Chern class of the pre-quantum line bundle.
The right input for the quantization procedure is
the group $DC_2^2$ of differential characters introduced by
Cheeger-Simons \cite{CheegerSimons1985}
(the indices of the notation $DC^2_2$ will be explained later). 
There are several definitions
of differential characters. The simplest one says that a
differential character of degree 2 (1 in Cheeger-Simons 
grading) is a pair $(\omega , \chi )$, where $\omega$ is a differential 
2-form and $\chi : Z_1 \rightarrow \mathbb{R}/\mathbb{Z}$ is a character 
of the group of smooth singular 1-cycles. This pair should satisfy the
following condition: 
\begin{equation}\nonumber
\chi ( \partial S )  \equiv \int_{S} \omega 
\quad \mathrm{mod} \ \mathbb{Z} \ ,
\end{equation}
for any smooth singular 2-chain $S$. One should think of $\omega$
as of a symplectic form (provided it is non-degenerate).
We use another definition of $DC_2^2$ which identifies it 
with the second cohomology group of
a certain complex $DC_2^{\bullet}$ involving both differential forms and
singular cochains (cf. \ref{DCChains}). In any case the 
crucial fact about differential characters is that they
classify isomorphism classes of Hermitian line bundles 
with connections on a manifold
$M$. Namely given such a bundle we can put $\omega$ to be the
curvature of the connection and $\chi$ to be its holonomy.
This map is a bijection and we call the inverse map
\begin{equation}\nonumber
DC_2^2 (M) \rightarrow \{\text{iso classes of Hermitian line bundles 
with connections on $M$} \}
\end{equation}
the prequantization. The actual quntization involves a choice of
polarization (for example a complex structure).
Our goal is to understand in
which sense the prequantization commutes with reduction.
The present paper provides the first step in this direction by explaining
prequantization in the equivariant setting (i.e. on stacks).
The reduction procedure is addressed in \cite{LermanMalkinReduction}.

Before proceeding with details of our construction of
equivariant prequantization, we would like (for technical reasons)
to switch from the category of Hermitian line bundles to the equivalent
category of principal $S^1$-bundles.

\subsection{Prequantization as a functor.}
Let $\Gamma_0 \leftleftarrows \Gamma_1$ be a Lie groupoid (one can
think, for example, of an action groupoid $M \leftleftarrows G \times M$).
An (equivariant) bundle on such a groupoid is a bundle on $\Gamma_0$
together with an isomorphism between the two pull-backs of the 
bundle to $\Gamma_1$. This definition leads us to an observation that
in order to understand equivariant prequantization
we have to include morphisms (in particular isomorphisms) between bundles
into our construction. So we promote the prequantization
from a bijection 
\begin{equation}\nonumber
DC_2^2 (M) \rightarrow \{\text{iso classes of Hermitian line bundles 
with connections on $M$} \}
\end{equation}
to a functor (equivalence of categories)
\begin{equation}\nonumber
\mathrm{Preq}_M : \mathcal{DC}_2^2 (M) \to \mathcal{DBS}^1 (M) \ ,
\end{equation}
where $\mathcal{DBS}^1 (M)$ is the category of principal 
$S^1$ bundles with connections over a manifold $M$ and
$\mathcal{DC}_2^2 (M)$ is a category which has $DC_2^2 (M)$
(i.e. second cohomology group of the complex $DC_2^{\bullet} (M)$)
as the set of isomorphism classes of objects. 
It is easy to see what the category $\mathcal{DC}_2^2 (M)$
should be: its objects are closed 2-cochains in $DC_2^{\bullet} (M)$
and the morphisms are 1-chains (see Section \ref{ChainStacks}
for a precise definition). This is an example of
a chain category, i.e. a category built from a complex 
of abelian groups. 

Observe that both categories $\mathcal{DC}_2^2 (M)$ and $\mathcal{DBS}^1 (M)$
can be pulled back under smooth maps of manifolds (and in fact,
$\mathrm{Preq}$ will intertwine the pull-back functors). Moreover
(and this is crucial for our argument) either of these categories
can be glued from its restrictions to open sets in a
covering of $M$. In other words there exists a stack
$\mathcal{DC}_2^2$ over the category of manifolds with the fiber 
over $M$ being $\mathcal{DC}_2^2 (M)$, and similar statement holds 
for $\mathcal{DBS}^1 (M)$. Once we know that 
$\mathcal{DBS}^1$ and $\mathcal{DC}_2^2$ are stacks, the general
abstract nonsense allows us to upgrade the functor
$\mathrm{Preq}$ from manifolds to orbifolds, Lie groupoids, 
and in fact arbitrary stacks (see Theorem \ref{PreQTheorem}). 
This solves the equivariant prequantization problem.

We should confess that we lifted the idea of applying chain categories
to classification problems from Hopkins-Singer paper
\cite{HopkinsSinger2005}. However they don't describe descent
properties of the categories involved and hence we could not quite 
fill-in details of their proofs. 
One can consider the present paper as a set of
exercises on some ideas of \cite{HopkinsSinger2005}.

\subsection{Chain stacks.}
The stack $\mathcal{DC}_2^2$ introduced above is an example of an 
Eilenberg-MacLane stack (cf. \cite{Toen2006}), 
which means that isomorphism classes of objects of its fiber 
$\mathcal{DC}_2^2 (M)$ over a manifold $M$ are in bijection with
some kind of cohomology (in this case differential characters). 
More familiar examples would be the stack
$\mathcal{H}^1 (C^{\bullet})$, which computes the first cohomology
of the singular cochain complex $C^{\bullet}$, or 
$\mathcal{H}^1 (\Omega^{\bullet})$ computing the first de Rham cohomology.
The importance of Eilenberg-MacLane stacks is that they allow one
to define equivariant cohomolgy, or more generally, stack cohomology,
in an intristic fashion
(without reference to an atlas). For example, the first
singular cohomology of an arbitrary stack $\mathcal{W}$
is defined to be the group of isomorphism classes of
objects of the functor category $\mathcal{H}om (\mathcal{W}, 
\mathcal{H}^1 (C^{\bullet}))$. This definition produces correct
answers for cohomology of manifolds and
equivariant cohomology of Lie groupoids. 

We would like to stress an important difference between 
what we call chain stacks and the general Eilenberg-MacLane stacks
defined in abstract algebraic topology. A chain stack is a 
presheaf of categories explicitly constructed from a complex of presheaves.
Being a stack is a condition on this presheaf of categories (and so
on the original complex of presheaves). An abstract definition of
an Eilenberg-MacLane stack involves stackification of a given presheaf
of categories. Roughly, working with chain stacks is similar to
working only with complexes of acyclic sheaves in sheaf cohomology.

Chain stacks provide a natural setup
for proving classification theorems. By a classification theorem
we mean a one-to-one correspondence between isomorphism 
classes of some kind of geometric objects on a manifold 
(for example, bundles with or without connections, gerbes, etc.) and 
some cohomology group of the manifold. Usually such theorems are
proved using \v{C}ech-type argument and proofs become increasingly
messy once one wants to deal with multiple covers or a group action.
A better approach is to show that the stack of objects we
want to classify is equivalent to a chain stack.
Since stacks are local, one just has to consider the case
of the manifold being an Euclidean space. Of course, the hard part
is to invent the cohomology theory represented by a chain stack
(see remarks after Theorem \ref{WeilTheorem}). 

An observant reader would notice that cohomology
groups in the above examples are all of low degrees. In fact higher 
cohomology is represented by higher stacks 
(cf. \citelist{\cite{Toen2006} \cite{Lurie2006}}). 
A higher stack is a presheaf of higher categories, or more precisely,
of simplicial sets, satisfying some descent condition.
The results of the present paper generalize to higher gerbs and
higher degree differential characters 
(cf. \cite{LermanMalkinNStacks}). In fact the proofs become
much cleaner when done in the language of simplicial sheaves.
However we want the present paper to be accessible to wider audience
and so we keep the exposition as elementary as possible
(at the expense of reproving many things which are well-known
and/or obvious in abstract algebraic topology/geometry)

\subsection{Past classification results.} 
Classification of principal $S^1$-bundles with or without connection
has a long history going back to work of 
Weil \cite{Weil1952}, who showed that the first Chern class
determines a line bundle up to an isomorphism. Kostant 
\cite{Kostant1970} described the prequantization map
from differential 2-forms with integer periods to 
isomorphism classes of line bundles 
with Hermitian connections. He also quantified the failure 
of the prequantization map to be a bijection. 
There are two closely related modifications of
prequantization which make it into a bijection:
Deligne cohomology (unpublished, see \cite{Beilinson1984}) 
and Cheeger-Simons differential characters \cite{CheegerSimons1985}.
We use the latter in the form modified by Hopkins-Singer
\cite{HopkinsSinger2005}.

Recently there has been an effort to generalize classification
theorems to equivariant setting. In the case of a group action
Brylinski \cite{Brylinski2000} proved Weil-type theorem using
good covers of the action groupoid, and
then this approach was extended to Deligne cohomology by
Gomi \cite{Gomi2005} (see also Lupercio-Uribe paper
\cite{LupercioUribe2006} for
the case of a finite group action).
Prequantization in the sense of Kostant
was generalized to arbitrary Lie groupoids by Behrend and Xu
\cite{BehrendXu2006}. Note that 
in the groupoid case the prequantization map
fails to be a bijection in a more complicated way that 
in the case of manifolds. 

The main results of the present paper are Weil and prequantization
theorems (\ref{ChernClassTheorem} and \ref{PreQTheorem}, resp.) 
for general stacks (over $\mathcal{M}an$). 
Our approach provides a new point of view on classification
theorems, which is useful for two reasons:
(1) the theory works for arbitrary stacks and is atlas-independent;
(2) it provides a bijection between isomorphism classes of bundles 
(with or without connections) and cohomology groups of certain
chain complexes (as opposed to sheaf cohomology groups).

\subsection{Warning: connections on equivariant bundles.}
\label{WarningSub}
One should be careful comparing various versions of equivariant 
classification/prequantization theorems. For example, in the case 
of a Lie group $G$ acting on a manifold $M$ one can consider equivariant
bundles with arbitrary, $G$-invariant, or $G$-basic connections.
Recall that a connection is $G$-basic iff it is $G$-invariant
and vanishes on vector fields generating the action of the
Lie algebra of $G$ on the total space of the bundle. In the present
paper we only consider basic connections on our bundles. These are
the connections which descend onto the quotient stack $[M/G]$,
make ``stacky'' sense, are atlas independent, etc.. The cases
of $G$-invariant and arbitrary connections on equivariant principal
$S^1$-bundles are treated in \cite{LermanMalkinReduction}.

\subsection{Structure of the paper.}
As mentioned above we wanted to keep the paper elementary as
possible. In particular, we reproduce many standard proofs from singular
homology theory and theory of stacks. Our goal is to define
cohomology stacks and to prove classification theorems
for principal $S^1$-bundles with or without connections.

Section \ref{BundleStacks} contains a review of stacks with
bundles being the main example.

Section \ref{ChainStacks} defines cohomology stacks for
singular and de Rham cohomology, and differential characters.

Section \ref{ChernFunctor} contains classification theorems
for principal $S^1$-bundles. 

\subsection*{Acknowledgements.} We would like to thank 
Matthew Ando and Thomas Nevins for valuable discussions. We also
thank the referee for careful reading of the paper.

\section{Stacks of bundles.}
\label{BundleStacks}

This section contains a brief review of the basic theory of
(differentiable) stacks. As examples we consider the stack of $S^1$-bundles
and the stack of $S^1$-bundles with connections. We refer the
reader to 
\citelist{\cite{BehrendXu2006} \cite{Metzler2003} \cite{Moerdijk2002}}
for a detailed exposition. Laumon and Moret-Bailly's book 
\cite{LaumonMoretBailly2000} develops the theory of algebraic stacks 
translatable into the language of differential geometry.

\subsection{Presheaves of groupoids.}
We denote by $\mathcal{M}an$ the category of differentiable 
manifolds (with smooth maps). 

Recall that a groupoid is a category such that any morphism is
invertible. A common notation for a groupoid is
$\Gamma_0 \leftleftarrows \Gamma_1$, where $\Gamma_0$ is the 
set of objects, $\Gamma_1$ is the set of morphisms, and the 
two arrows represent source and target maps. 
A groupoid is differentiable (or a Lie groupoid) if 
$\Gamma_0$ and $\Gamma_1$ are manifolds and all structure maps
(source, target, composition, inverse, identity) are smooth.

A (lax) presheaf of groupoids $\mathcal{X}$
over $\mathcal{M}an$ is a lax contravariant 2-functor from $\mathcal{M}an$
to the 2-category of groupoids. This means we have a groupoid
$\mathcal{X} (M)$ for each manifold $M$, a functor 
$f^* : \mathcal{X} (N) \rightarrow \mathcal{X} (M)$ 
for each smooth map $f: M \rightarrow N$, and
a natural transformation $f^* \circ g^* \cong (g \circ f)^*$
for each pair of composable smooth maps $f$ and $g$.

For example, given a manifold $M$, one can define a presheaf of groupoids
$[M]$ as follows. The objects of $[M]$ over a manifold
$N$ are smooth maps $N \rightarrow M$ and the only morphisms are the
identity ones. The pull-back functors are the usual pull-backs of
smooth maps.

To avoid (lax) 2-functors one can, instead
of presheaves of groupoids, consider an equivalent notion of
categories fibered in groupoids over $\mathcal{M}an$. Given such a 
category, its fibers form a presheaf of groupoids, and conversely,
one can build the total category from its fibers. 
See \cite{Metzler2003} for details.

Consider presheaves $\mathcal{X}$ and $\mathcal{Y}$ of groupoids
over $\mathcal{M}an$. The category 
$\mathcal{H}om (\mathcal{X} , \mathcal{Y})$ of morphisms
is defined as follows: the objects are maps 
$\mathcal{X} \rightarrow \mathcal{Y}$ of presheaves and the morphisms are 
natural transformations. For example, 2-Yoneda Lemma says that, 
given a manifold $M$ and an arbitrary presheaf of groupoids $\mathcal{X}$,
the category $\mathcal{H}om ([M] , \mathcal{X})$ is equivalent to
$\mathcal{X} (M)$. The equivalence functor is given (on objects) by:
\begin{equation}\nonumber
\mathrm{Ob} ( \mathcal{H}om ([M] , \mathcal{X}) ) \ni \quad f \mapsto 
f ( M \xrightarrow{\textrm{id}} M ) \quad 
\in \mathrm{Ob} (\mathcal{X} (M)) \ .
\end{equation}
We'll see later generalizations of this equivalence with $M$ replaced 
with a Lie groupoid or a stack.

In the present paper (presheaves of) groupoids play two roles. On one hand,
geometric objects we are interested in, such as base and the total space of 
a bundle, covering of a manifold, etc., will be thought of as (Lie) groupoids,
or (via 2-Yoneda Lemma) as presheaves of groupoids over $\mathcal{M}an$.
On the other hand, categories of bundles, chains of various complexes, etc.,
naturally form presheaves of groupoids over $\mathcal{M}an$. 

Some presheaves of groupoids are local (satisfy descent property) 
with respect to a topology on $\mathcal{M}an$ -- such presheaves 
are called stacks 
(see below for the precise definition). We use two pretopologies 
(in Grothendieck's sense). One is given by
open embeddings, the other by submersions. Since a submersion
has local sections, these pretopologies define the same topology. 

\subsection{Example: $\mathcal{BS}^1$.}
A typical example of a presheaf of groupoids is provided by
principal $S^1$-bundles. Namely, one defines a
presheaf $\mathcal{BS}^1$ of groupoids as follows:
\begin{itemize}
\item
Given a manifold $M$, objects of $\mathcal{BS}^1 (M)$ are principal 
$S^1$-bundles $P \rightarrow M$. Abusing notation we often specify
only the total space $P$. 
\item
A morphism in $\mathcal{BS}^1 (M)$ from $P$ to $P'$ is a smooth map
$\phi: P \rightarrow P'$ such that the following diagram commutes:
\begin{equation}\nonumber
\xymatrix{
P \ar[r]^{\phi} \ar[d] & P' \ar[d] \\
M \ar[r]^{\mathrm{id}} & M
}
\end{equation}
We call such $\phi$ a bundle map.
\item
Given a smooth map $f : M \rightarrow N$ the functor 
$f^* : \mathcal{BS}^1 (N) \rightarrow \mathcal{BS}^1 (M)$
is the pull-back of $S^1$-bundles (note that pull-backs
are only defined up to isomorphism).
\end{itemize}
Note that 2-Yoneda Lemma implies that the category $\mathcal{BS}^1 (M)$ of
principal $S^1$-bundles on a manifold $M$ is equivalent to 
the morphisms category 
$\mathcal{H}om ([M], \mathcal{BS}^1)$. This is to be compared with the 
definition of the classifying space $BS^1$ in topology: $BS^1$ is a 
topological space such that the homotopy classes of maps $M \rightarrow BS^1$ 
are in bijection with equivalence classes of $S^1$-bundles on $M$.

\subsection{Equivariant objects of a presheaf.}
\label{EquivariantPresheaves}
Consider a Lie groupoid 
\mbox{$\Gamma_{\bullet} = \{ \Gamma_0 \leftleftarrows \Gamma_1 \}$}. 
The nerve of 
$\Gamma_{\bullet}$ is the simplicial manifold
\begin{equation}\nonumber
\xymatrix{
{\Gamma_0} &
\ar@<0.3ex>[l]^{\partial_0}
\ar@<-0.7ex>[l]_{\partial_1}  
{\Gamma_1} &
\ar@<0.9ex>[l]^{\partial_0}
\ar@<-0.1ex>[l]
\ar@<-1.1ex>[l]_{\partial_2}  
{\Gamma_2} &
\ar@<1.3ex>[l]^{\partial_0}
\ar@<0.3ex>[l]
\ar@<-0.7ex>[l]
\ar@<-1.7ex>[l]_{\partial_3}  
{\cdots}}
\end{equation}
where
\begin{equation}\nonumber
\Gamma_n = \underbrace{
\Gamma_1 \times_{\Gamma_0} \Gamma_1 \times_{\Gamma_0}
\ldots \times_{\Gamma_0} \Gamma_1}_n
\end{equation}
and arrows are the canonical projections. We omit face maps (in particular,
the identity map for the groupoid) in the above diagram. We also don't  
distinguish between a groupoid and the associated simplicial manifold and denote
both by $\Gamma_{\bullet}$.

Now let $\mathcal{X}$ be a presheaf of groupoids over $\mathcal{M}an$.
We define a category $\mathcal{X} (\Gamma_\bullet )$ of
objects of $\mathcal{X}$ on $\Gamma_\bullet$ 
(or $\Gamma_\bullet$-equivariant objects of $\mathcal{X}$)
as follows:
\begin{itemize}
\item
An object of $\mathcal{X} (\Gamma_\bullet )$ is a pair $(x, \phi )$,
where $x$ is an object of $\mathcal{X} (\Gamma_0 )$ and
$\phi : \partial_0^* x \rightarrow \partial_1^* x$ is an 
isomorphism in $\mathcal{X} (\Gamma_1 )$ satisfying a cocycle 
condition on $\mathcal{X} (\Gamma_2 )$.
\item
A morphism from $(x, \phi )$ to $(x', \phi' )$ is a morphism
$\xi : x \rightarrow x'$ in $\mathcal{X} (\Gamma_0 )$ such that
$\phi' \circ \partial_0^* \xi = \partial_1^* \xi \circ \phi$
in $\mathcal{X} (\Gamma_1 )$.
\end{itemize}
If the base groupoid is just a manifold
$M_{\bullet} = \xymatrix@1{
M &\ar@<0.2ex>[l]^{\mathrm{id}} \ar@<-0.8ex>[l]_{\mathrm{id}}  M}$
then the category $\mathcal{X} ( M_{\bullet} ) $ is equivalent to the
category $\mathcal{X} (M)$.
Another important example is provided by an action
of a Lie group $G$ on a manifold $M$. The associated (action) 
groupoid is 
$\xymatrix@1{{M} &\ar@<0.2ex>[l]^{p} \ar@<-0.8ex>[l]_{a} G \times M}$,
where the arrows $a$ and $p$ represent the action and 
the canonical projection respectively. 
The category $\mathcal{X} ( M \leftleftarrows G \times M )$
is the category of $G$-equivariant objects of $\mathcal{X} (M)$. 
For example $\mathcal{BS}^1 ( M \leftleftarrows G \times M )$ is
the category of $G$-equivariant $S^1$-bundles on $M$.

\subsection{Descent.} 
We turn our attention to the topological structure of the
site $\mathcal{M}an$. A presheaf of groupoids is a stack if is is 
local; in other words, if global objects are glued from local ones.
To make it precise one uses equivariant objects for a covering
groupoid.

Given a covering $U \rightarrow M $ of a manifold $M$, we 
denote by $U_{\bullet}$ the Lie groupoid 
$U \times_M U \rightrightarrows U$. The groupoid
$U_{\bullet}$ is called the covering groupoid (corresponding to
the covering $U \rightarrow M$). 
In the open subsets topology we have 
$U = \bigsqcup_{\alpha} U_{\alpha}$,
where $\{ U_{\alpha} \}_{\alpha}$ is a covering of $M$
by open subsets, and then
$U \times_M U = \bigsqcup_{\alpha , \beta} U_{\alpha} \cap U_{\beta}$.

If $\mathcal{X}$ is a presheaf of groupoids and $U_{\bullet}$ is
a covering groupoid, then the groupoid $\mathcal{X} ( U_{\bullet} )$
is called the descent data (of $\mathcal{X}$ with respect to
the covering $U \rightarrow M$). The descent data is effective if
it is equivalent to the groupoid $\mathcal{X} (M)$ (more precisely,
if the natural restriction functor $\mathcal{X} (M) \rightarrow 
\mathcal{X} ( U_{\bullet} )$ is an equivalence of categories).
A presheaf of groupoids is a stack if its descent data is effective on any
manifold with respect to any covering (it is enough to consider 
coverings by open subsets).

\begin{proposition}
The presheaf $\mathcal{BS}^1$ is a stack. 
\end{proposition}
\begin{proof} 
Given a covering $M = \bigcup U_{\alpha}$ of $M$ by open subsets,
an object of the descent data is a collection of bundles on
the subsets $U_{\alpha}$ together with bundle isomorphisms on 
the intersections $U_{\alpha} \cap U_{\beta}$. These isomorphisms 
satisfy cocycle conditions on triple
intersections and therefore define an equivalence
relation on the disjoint union of the total spaces. 
The equivalence classes of this relation form the total space of
an $S^1$-bundle on $M$. It is easy
to see that the gluing procedure is a functor
quasi-inverse (inverse up to natural transformations) to restriction
functor.
\end{proof}

\subsection{Differentiable stacks.}
In the last several subsections we considered objects of a stack
on manifolds and, more generally, on Lie groupoids. Now we want to explain
what an object of a stack (say, an $S^1$-bundle) on another
stack is. We start with recalling in what sense stacks are generalizations
of Lie groupoids.

Given a Lie groupoid $\Gamma_0 \leftleftarrows \Gamma_1$ 
(for example, associated to a Lie group action 
$M \leftleftarrows G \times M$), one can define the quotient
stack $[\Gamma_0 / \Gamma_1 ]$ as the classifying stack of 
principal $\Gamma_{\bullet}$-bundles. Recall that a principal 
$\Gamma_{\bullet}$-bundle
on a manifold $M$ is a surjective submersion $\pi: P \rightarrow M$ 
together with a right action of $\Gamma_{\bullet}$ on $P$ which commutes with
the projection $\pi$ and is free and transitive on fibers. 
A right action of $\Gamma_{\bullet}$ on $P$ is given by an anchor map 
$P \rightarrow \Gamma_0$ and an action map 
$P \times_{\Gamma_0} \Gamma_1 \rightarrow P$ 
satisfying a set of axioms. We write the action map as
$(p,g) \mapsto pg$. Since the action commutes with $\pi$
(i.e. $\pi (pg)=\pi (p)$), the assignment 
$(p,g) \mapsto (p, pg)$ defines a map
$P \times_{\Gamma_0} \Gamma_1 \rightarrow P \times_M P$ .
The action is free and transitive on fibers of $\pi$ iff this map is
a diffeomorphism. 

The above discussion can be summarized by saying 
that a principal $\Gamma_{\bullet}$-bundle over $M$ is a commutative diagram:
\begin{equation}\label{GammaBundle}
\xymatrix{
{M} & 
\ar[l] \ar[d] {P} &
\ar@<0.3ex>[l] \ar@<-0.7ex>[l] \ar[d] {P \times_M P} &
\ar@<0.9ex>[l] \ar@<-0.1ex>[l] \ar@<-1.1ex>[l] \ar[d]
{P \times_M P \times_M P} &
\ar@<1.3ex>[l]
\ar@<0.3ex>[l]
\ar@<-0.7ex>[l]
\ar@<-1.7ex>[l]  
{\cdots} \\
& {\Gamma_0} &
\ar@<0.3ex>[l]
\ar@<-0.7ex>[l]  
{\Gamma_1} &
\ar@<0.9ex>[l]
\ar@<-0.1ex>[l]
\ar@<-1.1ex>[l]  
{\Gamma_2} &
\ar@<1.3ex>[l]
\ar@<0.3ex>[l]
\ar@<-0.7ex>[l]
\ar@<-1.7ex>[l]  
{\cdots}}
\end{equation}
Note that the top row is a covering 
groupoid for $M$ (in the submersion topology). 
Also, locally on $M$ there exists a section 
$M \rightarrow \Gamma_0$ trivializing the bundle.
One easily defines morphisms of $\Gamma_{\bullet}$-bundles
and the pull-back functor with respect to smooth maps.
Hence $\Gamma_{\bullet}$-bundles form a 
presheaf  $[\Gamma_0 / \Gamma_1 ]$ of groupoids, 
and it is easy to check that  $[\Gamma_0 / \Gamma_1 ]$ is a stack.
For example, $\mathcal{BS}^1 = [\mathrm{pt}/S^1]$.

A stack $\mathcal{W}$ is a differentiable if it is equivalent to
a stack $[ \Gamma_0 / \Gamma_1 ]$ for some groupoid 
$\Gamma_{\bullet}$. Equivalently, $\mathcal{W}$ is differentiable
if it has a smooth atlas (surjective representable morphism)
$[\Gamma_0] \rightarrow \mathcal{W}$. 
Recall (cf. \citelist{\cite{Heinloth2005} \cite{Metzler2003}})
that a morphism $[X] \rightarrow \mathcal{V}$ from a manifold
to a stack is surjective representable if for any morphism 
$[Y] \rightarrow \mathcal{V}$ the stack $[Y] \times_{\mathcal{V}} [X]$ 
is representable by a non-empty manifold. 
Given an atlas $[\Gamma_0] \rightarrow \mathcal{W}$, one lets 
$\Gamma_1$ to be a representative of $[\Gamma_0] \times_{\mathcal{W}} [\Gamma_0]$
(i.e. $[\Gamma_1] \cong [\Gamma_0] \times_{\mathcal{W}} [\Gamma_0]$). 
Conversely, the trivial $\Gamma_{\bullet}$-bundle over
$\Gamma_0$ defines an atlas 
$[\Gamma_0] \rightarrow [ \Gamma_0 / \Gamma_1]$.
Abusing terminology we will call the whole simplicial
manifold $\Gamma_{\bullet}$ the atlas. 

The following standard proposition 
shows that morphisms of stacks are local (i.e. can be
described using atlases).

\begin{proposition}\label{XGHomGX}
If $\mathcal{X}$ is an arbitrary stack and $\mathcal{W}$
is a differentiable stack with an atlas $\Gamma_{\bullet}$,
then the following three categories: 
$\mathcal{H}om ( \mathcal{W},  \mathcal{X} )$, 
$\mathcal{H}om ( [ \Gamma_1 / \Gamma_0 ],  \mathcal{X} )$, and 
$\mathcal{X} (\Gamma_{\bullet})$, are equivalent to each other.
\end{proposition}
\begin{proof}
The first two categories are obviously equivalent (since 
$\mathcal{W}$ is equivalent to $[ \Gamma_1 / \Gamma_0 ]$).
It remains to be shown that they are equivalent to
$\mathcal{X} (\Gamma_{\bullet})$. We provide functors in both directions. 
First, precomposition with the atlas map 
$[\Gamma_0] \rightarrow \mathcal{W}$ 
defines a functor from 
$\mathcal{H}om ( \mathcal{W},  \mathcal{X} )$ to
$\mathcal{H}om ( [ \Gamma_0 ],  \mathcal{X} ) \stackrel{\text{2-Yoneda}}{\cong}
\mathcal{X} ( \Gamma_0 )$ together with a natural transformation
(satisfying a natural cocycle condition)
between the two pull-backs of this functor to  
$\mathcal{H}om ( [ \Gamma_0 ] \times_\mathcal{W} [ \Gamma_0 ],  
\mathcal{X} ) \cong
\mathcal{X} (\Gamma_1)$. These data is the same as a functor 
from $\mathcal{H}om ( \mathcal{W},  \mathcal{X} )$ to
$\mathcal{X} ( \Gamma_{\bullet} )$. Up to this point the argument works for any
presheaf of categories $\mathcal{X}$; however
existence of a quasi-inverse functor 
$\mathcal{X} ( \Gamma_{\bullet} ) \rightarrow
\mathcal{H}om ( [ \Gamma_0 / \Gamma_1 ],  \mathcal{X} )$
requires
the assumption that $\mathcal{X}$ is a stack. 
Let us describe this functor on objects (extension to morphisms
is clear). Given an object of 
$\mathcal{X} ( \Gamma_{\bullet} )$ and a 
$\Gamma_{\bullet}$-bundle over a manifold $M$, we
want to produce an object of $\mathcal{X} (M)$. So we pull back
$\mathcal{X} ( \Gamma_{\bullet} )$ to 
$\mathcal{X} (P_{\bullet} )$ along the vertical projection in
\eqref{GammaBundle}. The result is a descent data object 
of $\mathcal{X}$ for the covering groupoid $P_{\bullet}$ of $M$, 
which determines an object of $\mathcal{X} (M)$ since $\mathcal{X}$ 
is a stack.
\end{proof}

The above proposition motivates the following definition, which extends
the notion of an object of a stack on a manifold (or on a groupoid).
Given two stacks $\mathcal{X}$ and $\mathcal{W}$
we define the category $\mathcal{X} (\mathcal{W})$ of objects
of $\mathcal{X}$ on  $\mathcal{W}$ to be the category 
$\mathcal{H}om ( \mathcal{W}, \mathcal{X} )$. Note that up to an
equivalence this category depends only on the equivalence
classes of $\mathcal{W}$ and $\mathcal{X}$. If $\mathcal{W}$
is differentiable with an atlas $\Gamma_{\bullet}$ then 
the category $\mathcal{X} (\mathcal{W})$ is equivalent to 
the category of $\Gamma_{\bullet}$-equivariant objects of 
$\mathcal{X}$. 
The new definition however is independent of the atlas
(up to an equivalence).

\subsection{Stackification and $\mathcal{BS}_{\mathrm{triv}}^1$.}

Given a presheaf of groupoids $\mathcal{X}_0$ over $\mathcal{M}an$,
there exists unique (up to an equivalence) stack 
$\mathcal{X}$  together
with a functor $s: \mathcal{X}_0 \rightarrow \mathcal{X}$ 
satisfying the following universal property: 
precomposition with $s$ is an equivalence of categories
$\mathcal{H}om (\mathcal{X}, \mathcal{Y}) \cong 
\mathcal{H}om (\mathcal{X}_0, \mathcal{Y})$
for any stack $\mathcal{Y}$. One can construct a
stackification as the union of the descent data of
$\mathcal{X}_0$ with respect to all coverings up to equivalence
corresponding to refinements of coverings.

As an example consider a full sub-presheaf $\mathcal{BS}_{\mathrm{triv}}^1$ of 
$\mathcal{BS}^1$ consisting of trivial bundles, i.e. for a 
manifold $M$ the only object of $\mathcal{BS}_{\mathrm{triv}}^1 (M)$ is the
trivial bundle $M \times S^1$.
Recall that a sub-presheaf of groupoids is called full if the set of morphisms
between any two objects is the same as in the original presheaf.
In particular, the morphisms in $\mathcal{BS}_{\mathrm{triv}}^1 (M)$ 
are functions $M \rightarrow S^1$.
The presheaf $\mathcal{BS}_{\mathrm{triv}}^1$ is not a stack 
but we have the following
\begin{proposition}
The stackification of $\mathcal{BS}_{\mathrm{triv}}^1$ is $\mathcal{BS}^1$.
\end{proposition}
\begin{proof}
Any $S^1$-bundle is locally trivial.
\end{proof}

\subsection{The classifying stack $\mathcal{DBS}^1$.}
Now we want to equip our bundles with connections. 
So let $\mathcal{DBS}^1$ be
the following presheaf of groupoids over $\mathcal{M}an$:
\begin{itemize}
\item
Given a manifold $M$, objects of $\mathcal{DBS}^1 (M)$ are principal 
$S^1$-bundles with connections over $M$, i.e. pairs $(P,A)$, where
$P \rightarrow M$ is a bundle and $A \in \Omega^1 (P)$ is a connection.
\item
Given two objects $(P, A)$ and $(P' , A')$ of $\mathcal{DBS}^1 (M)$,
a morphism from $(P, A)$ to $(P' , A')$ is a bundle map 
$\phi: P \rightarrow P'$ such that $\phi^* A' =A$.
\item
Given a smooth map $f : M \rightarrow N$ the functor 
$f^* : \mathcal{DBS}^1 (N) \rightarrow \mathcal{DBS}^1 (M)$
is the pull-back of bundles with connections.
\end{itemize}
As with any presheaf of groupoids we can consider
categories of equivariant objects of $\mathcal{DBS}^1$. 
For example, in the case of an action groupoid, the category
\mbox{$\mathcal{DBS}^1 ( M \leftleftarrows G \times M )$} is
the category of $G$-equivariant principal $S^1$-bundles on $M$ with basic
connections (cf. subsection \ref{WarningSub}).

\begin{proposition}
The presheaf $\mathcal{DBS}^1$ is a stack.
\end{proposition}
\begin{proof}
Similar to the case of $\mathcal{BS}^1$. Note that gluing isomorphisms 
preserve connections. Hence connections on local bundles
determine a connection on the global bundle. 
\end{proof}

\subsection{The presheaf $\mathcal{DBS}_{\mathrm{triv}}^1$.}
Similar to $\mathcal{BS}_{\mathrm{triv}}^1$, we introduce a full presheaf 
$\mathcal{DBS}_{\mathrm{triv}}^1$ of $\mathcal{DBS}^1$
consisting of trivial bundles with connections, i.e. for a 
manifold $M$ objects of $\mathcal{DBS}_{\mathrm{triv}}^1 (M)$ are pairs 
$(M \times S^1, \alpha + d \theta)$, where $\alpha \in \Omega^1 (M)$.
The morphisms in $\mathcal{DBS}_{\mathrm{triv}}^1$ are 
locally constant functions $M \rightarrow S^1$.

\begin{proposition}
The stackification of $\mathcal{DBS}_{\mathrm{triv}}^1$ is $\mathcal{DBS}^1$.
\end{proposition}
\begin{proof}
Any $S^1$-bundle is locally trivial and any connection on a trivial
bundle is of the form $\alpha+d \theta$.
\end{proof}

\section{Chain stacks.}
\label{ChainStacks}

In this section we introduce chain stacks associated to complexes
of presheaves. In particular we define differential characters.

\subsection{Categories built from complexes.}
\label{ChainCategories}

Let $A^{\bullet} = \{ A^\bullet \xrightarrow{d} A^{\bullet+1} \}$ 
be a complex of abelian groups. 
All complexes appearing in this paper are assumed
to have $A^n = 0$ for $n<0$. We fix an integer $n \geq 0$ and 
define a category $\mathcal{H}^n ( A^\bullet )$ as follows:
\begin{itemize}
\item
Objects are $n$-cocycles: $z \in A^n$, $dz = 0$.
\item
Morphisms are (n-1)-cochains up to (n-2)-cochains. The set of morphisms
from $z$ to $z'$ is   
\begin{equation}\nonumber
\{ b \in A^{n-1} \ | \ db = z'-z \} / (b \sim b +dc , \ c \in A^{n-2}) \ .
\end{equation}
Composition of morphisms is the addition in $A^{n-1}$.
\end{itemize}
The map $\mathcal{H}^n$ can be extended to a functor from the category
of complexes (and chain maps) to the 2-category of categories:
given a map $f : A^\bullet \rightarrow B^\bullet$ we have a functor
$\mathcal{H}^n (f) : \mathcal{H}^n ( A^\bullet )
\rightarrow \mathcal{H}^n ( B^\bullet )$ defined as follows:
\begin{itemize}
\item
On objects: $\mathcal{H}^n (f) (z) = f(z)$ for $z \in A^n$, $dz=0$.
\item
On morphisms: $\mathcal{H}^n (f) ([b]) = [f(b)]$ for $[b] 
\in A^{n-1}/d A^{n-2}$.
\end{itemize}
Moreover, let $k : A^\bullet \rightarrow B^{\bullet-1}$ 
be a chain homotopy between two chain maps
$f$, $g : A^\bullet \rightarrow B^\bullet$
(i.e. $g-f=dk+kd$). Then we have a natural transformation 
$\mathcal{H}^n (k) : \mathcal{H}^n ( f )
\rightarrow \mathcal{H}^n ( g )$ given by
$\bigl( \mathcal{H}^n (k) \bigr) (z) = [k(z)]$ for $z \in A^n$, $dz=0$.

Let us make a few remarks about this construction:
\begin{itemize}
\item
The category $\mathcal{H}^n ( A^\bullet )$ is a groupoid.
\item
The set of objects of $\mathcal{H}^n ( A^\bullet )$ is an abelian group. 
\item
The set of isomorphism classes of objects of $\mathcal{H}^n ( A^\bullet )$ 
is the $n^{\mathrm{th}}$ cohomology group $H^n ( A^\bullet )$ of 
the complex $A^\bullet$.
\item
The automorphism group of any object of $\mathcal{H}^n ( A^\bullet )$ 
is $H^{n-1} ( A^\bullet )$.
\item
The category $\mathcal{H}^0 ( A^\bullet )$ is discrete (the only 
morphisms are identity ones): we assume our complexes to
be trivial in negative degrees.
\end{itemize}
It follows that, if $f : A^\bullet \rightarrow B^{\bullet}$ is a
quasi-isomorphism (i.e. induces isomorphism in cohomology), 
then $\mathcal{H}^n(f)$ is an equivalence
of categories. The reason is that a functor from a groupoid to a groupoid
is an equivalence if it preserves isomorphism classes of objects and
their automorphism groups.

Now let $F^{\bullet}$ be a complex of presheaves of abelian groups over
$\mathcal{M}an$. This means we have a complex 
$F^\bullet (M)$ of abelian groups for 
each manifold $M$ and a pull-back map 
$F^\bullet (N) \rightarrow F^\bullet (M)$ for each smooth
map $M \rightarrow N$.
We can apply the functor $\mathcal{H}$ to
$F^\bullet$ and obtain a presheaf $\mathcal{H} (F^\bullet)$ of
groupoids. 

\subsection{Equivariant cochains.}
Given a Lie groupoid $\Gamma_{\bullet}$ and a complex $F^{\bullet}$ 
of pre-sheaves of abelian groups we consider the double complex 
$F^{\bullet} ( \Gamma_{\bullet} )$ with the differentials 
$d : F^{\bullet} ( \Gamma_{\bullet} ) \rightarrow 
F^{\bullet + 1} ( \Gamma_{\bullet} )$ and 
$\delta : F^{\bullet} ( \Gamma_{\bullet} ) \rightarrow 
F^{\bullet} ( \Gamma_{\bullet +1} )$:
\begin{equation}\nonumber
\xymatrix{
{\vdots} & & {\vdots} & {\vdots} & {\vdots} & \\
{F^2} \ar[u]^{d} & &
{F^2 ( \Gamma_0 )} 
\ar[r]^{\delta} \ar[u]^{d} &
{F^2 ( \Gamma_1 )} 
\ar[r]^{\delta} \ar[u]^{-d} &
{F^2 ( \Gamma_2 )} 
\ar[r]^{\delta} \ar[u]^{d} &
{\cdots}
\\
{F^1} \ar[u]^{d} & &
{F^1 ( \Gamma_0 )} 
\ar[r]^{\delta} \ar[u]^{d} &
{F^1 ( \Gamma_1 )} 
\ar[r]^{\delta} \ar[u]^{-d} &
{F^1 ( \Gamma_2 )} 
\ar[r]^{\delta} \ar[u]^{d} &
{\cdots}
\\
{F^0} \ar[u]^{d} & &
{F^0 ( \Gamma_0 )} 
\ar[r]^{\delta} \ar[u]^{d} &
{F^0 ( \Gamma_1 )} 
\ar[r]^{\delta} \ar[u]^{-d} &
{F^0 ( \Gamma_2 )} 
\ar[r]^{\delta} \ar[u]^{d} &
{\cdots}
\\ \\ & &
{\Gamma_0} &
\ar@<0.3ex>[l]^{\partial_0}
\ar@<-0.7ex>[l]_{\partial_1}  
{\Gamma_1} &
\ar@<0.9ex>[l]^{\partial_0}
\ar@<-0.1ex>[l]
\ar@<-1.1ex>[l]_{\partial_2}  
{\Gamma_2} &
\ar@<1.3ex>[l]^{\partial_0}
\ar@<0.3ex>[l]
\ar@<-0.7ex>[l]
\ar@<-1.7ex>[l]_{\partial_3}  
{\cdots}
}
\end{equation}
The differential $d$ is given by the original differential in the
complex of presheaves $F^{\bullet}$ and the differential
$\delta$ is constructed using the structure of a simplicial
manifold on the nerve of $\Gamma_{\bullet}$:
\begin{equation}\nonumber
\delta = \sum_i (-1)^i \partial_i^*
\end{equation}
The cohomology of the total complex 
$( F(\Gamma)_{\mathrm{tot}}^{\bullet} , d_{\mathrm{tot}} )$
associated to the above double complex is called 
the $\Gamma_{\bullet}$-equivariant cohomology of
$F^{\bullet}$ (cf. \cite{Behrend2004}). 

Presently we have two categories associated to a complex of presheaves
$F^{\bullet}$, a Lie groupoid $\Gamma_{\bullet}$, and an integer $n \geq 0$:
the category 
$\bigl( \mathcal{H}^n ( F^\bullet ) \bigr) ( \Gamma_{\bullet} )$ 
of $\Gamma_{\bullet}$-equivariant objects of the presheaf
$\mathcal{H}^n ( F^\bullet )$ 
(cf. \ref{EquivariantPresheaves}) and the chain category
of the total complex $\mathcal{H}^n ( F(\Gamma)_{\mathrm{tot}}^{\bullet} )$.
These categories are not equivalent for a general $n$,
but $n=0$ and $n=1$ are exceptional cases.
\begin{proposition}\label{H01Total}
The following categories are isomorphic:
\begin{equation}
\begin{split}
\bigl( \mathcal{H}^0 ( F^\bullet ) \bigr) ( \Gamma_{\bullet} ) = 
\mathcal{H}^0 ( F(\Gamma)_{\mathrm{tot}}^{\bullet} ) \\
\bigl( \mathcal{H}^1 ( F^\bullet ) \bigr) ( \Gamma_{\bullet} ) = 
\mathcal{H}^1 ( F(\Gamma)_{\mathrm{tot}}^{\bullet} ) 
\end{split}
\end{equation}
\end{proposition}
\begin{proof}
The isomorphisms follow directly from definitions. 
Here are the details.

Degree $0$ case. Both 
$(\mathcal{H}^0 ( F^\bullet ) \bigr) ( \Gamma_{\bullet} )$ and
$\mathcal{H}^0 ( F(\Gamma)_{\mathrm{tot}}^{\bullet} )$ are discrete 
categories with the set of objects 
$\{ z \in F^0 ( \Gamma_0 ) \ | \ dz=0 , 
\ \partial_0^* z  = \partial_1^* z \}$
and no morphisms except identity ones.

Degree $1$ case. Objects of 
$\mathcal{H}^1 ( F(\Gamma)_{\mathrm{tot}}^{\bullet} )$ are pairs
\mbox{$(z,t) \in F^1 ( \Gamma_0 ) \oplus F^0 ( \Gamma_1 )$} such that
\begin{align}
dz &= 0 \label{dz} \\ 
\partial_0^* z  - \partial_1^* z &= d t \label{dzdt} \\
\partial_0^* t  - \partial_1^* t + \partial_2^* t &= 0 \label{dt}
\end{align}
Now observe that \eqref{dz} means $z$ is an object of 
$\bigl( \mathcal{H}^1 ( F^{\bullet}) \bigr) ( \Gamma_0 )$,
\eqref{dzdt} means $t$ is a morphism from
$\partial_0^* z$ to $\partial_1^* z$ in
$\bigl( \mathcal{H}^1 ( F^{\bullet} ) \bigr) ( \Gamma_1 )$, and
\eqref{dt} is the cocycle condition for this
morphism. In other words, the pair $(z,t)$ is an object of 
$\bigl( \mathcal{H}^1 ( F^\bullet ) \bigr) ( \Gamma_{\bullet} )$.
A morphism from $(z,t)$ to $(z',t')$ in either 
$\mathcal{H}^1 ( F(\Gamma)_{\mathrm{tot}}^{\bullet} )$ or
$\bigl( \mathcal{H}^1 ( F^\bullet ) \bigr) ( \Gamma_{\bullet} )$ 
is an element $b \in F^0 ( \Gamma_0 )$ such that $db = z-z'  , 
\ \partial_0^* b  - \partial_1^* b  = t - t'$.
\end{proof}

\subsection{Example: singular cochains.}\label{SingularChains}
Given a manifold $M$, we consider the complex $C_{\bullet} (M)$ 
of smooth singular chains: the group
$C_n (M)$ is a free abelian group generated by all
smooth maps from the $n$-simplex to $M$, and the differential is 
the boundary operator. The dual complex 
is the complex of smooth singular cochains  $C^{\bullet} (M)$.
Cochains can be naturally pulled back with respect to smooth
maps, and we obtain a complex $C^{\bullet}$ of presheaves of abelian groups.

It is easy to see that $C^0$ and $H^0 ( C^{\bullet} )$ are sheaves.
On the other hand, $H^1 ( C^\bullet )$ is not a sheaf. In fact not even
$C^1$ is a sheaf (this is a well-known technical issue in singular 
(co)homology theory). However we have the following

\begin{proposition}[Mayer-Vietoris]
The presheaf of groupoids $\mathcal{H}^1 (C^\bullet)$ is a stack.
\end{proposition}

The presheaf of isomorphism classes of objects of 
$\mathcal{H}^1 (C^\bullet)$ is the first singular cohomology presheaf.
One can phrase the above proposition as ``first cohomology is local
if one thinks of it as a stack''.

\begin{proof}
We have to show that given an open covering 
$U=\bigsqcup_{\alpha} U_{\alpha} \xrightarrow{\upsilon} M$ 
of a manifold $M$ with the associated
(\v{C}ech) groupoid $U_{\bullet}$, the restriction map is an
equivalence between groupoids 
$\bigl( \mathcal{H}^1 ( C^{\bullet} ) \bigr) (M) =
\mathcal{H}^1 ( C^{\bullet} (M))$ and
$\bigl( \mathcal{H}^1 ( C^{\bullet} )\bigr) ( U_{\bullet} )$. 
By Proposition \ref{H01Total} we can replace the 
latter with $\mathcal{H}^1 ( C(U)_{\mathrm{tot}}^{\bullet} )$.
So it remains to be shown that the restriction map
$\upsilon^* : C^{\bullet} (M) \rightarrow C^{\bullet} (U)$ 
followed by the inclusion $\epsilon$ of the complex
$C^{\bullet} (U)$ into the double complex 
$C( U)^{\bullet}_{\mathrm{tot}}$ induces an equivalence
of categories $\mathcal{H}^1 (\epsilon \circ \upsilon^*) : 
\mathcal{H}^1 ( C^{\bullet} (M)) \rightarrow 
\mathcal{H}^1 ( C (U)_{\mathrm{tot}}^{\bullet} )$.
This follows from a version of the Mayer-Vietoris theorem
saying that $\epsilon \circ \upsilon^*$ is a quasi-isomorphism.

Let us explain details of the proof showing that 
the restriction $\mathcal{H}^1 (\epsilon \circ \upsilon^*) : 
\mathcal{H}^1 ( C^{\bullet} (M)) \rightarrow 
\mathcal{H}^1 ( C (U)_{\mathrm{tot}}^{\bullet} )$
is an equivalence.
The argument is standard in singular homology
theory but we want to rephrase it in terms of chain categories
and functors.

\emph{Step 1.} 
We consider a complex $C^{\bullet}(M, U)$ of ``small''
cochains, i.e. cochains that take values only on those simplexes
$\Delta \rightarrow M$ that can be lifted to $\Delta \rightarrow U$
(in other words such that the image of $\Delta$ is contained in
a single open set $U_{\alpha}$).
The standard subdivision argument (see \cite{Vick1994}) shows that
the restriction map $C^{\bullet}(M) \rightarrow C^{\bullet}(M, U)$ 
has an inverse up to homotopy, which
according to Subsection \ref{ChainCategories} corresponds to the
restriction functor $\mathcal{H}^1 ( C^{\bullet}(M) ) \rightarrow 
\mathcal{H}^1 ( C^{\bullet}(M, U) )$ being
equivalence of categories. 
More precisely, in \cite{Vick1994} a homotopy between
small and arbitrary chains is given by subdivision. We use
the dual homotopy between cochains. In any case it remains to be shown 
that the restriction $\mathcal{H}^1 ( C^{\bullet}(M, U) ) \rightarrow
\mathcal{H}^1 ( C(U)_{\mathrm{tot}}^{\bullet} )$ is an equivalence.

\emph{Step 2.} Observe that 
$\mathcal{H}^1 ( C(U)_{\mathrm{tot}}^{\bullet} )$
is the chain category of the total complex of the following
double complex
\begin{equation}\label{SingularDoubleComplex}
\xymatrix{
{\vdots} & {\vdots} & {\vdots} & \\
{C^2 ( U )} 
\ar[r]^-{\delta} \ar[u]^{d} &
{C^2 ( U \times_M U )} 
\ar[r]^-{\delta} \ar[u]^{-d} &
{C^2 ( U \times_M U \times_M U )} 
\ar[r]^-{\delta} \ar[u]^{d} &
{\cdots}
\\
{C^1 ( U )} 
\ar[r]^-{\delta} \ar[u]^{d} &
{C^1 ( U \times_M U )} 
\ar[r]^-{\delta} \ar[u]^{-d} &
{C^1 ( U \times_M U \times_M U )} 
\ar[r]^-{\delta} \ar[u]^{d} &
{\cdots}
\\
{C^0 ( U )} 
\ar[r]^-{\delta} \ar[u]^{d} &
{C^0 ( U \times_M U )} 
\ar[r]^-{\delta} \ar[u]^{-d} &
{C^0 ( U \times_M U \times_M U )} 
\ar[r]^-{\delta} \ar[u]^{d} &
{\cdots}}
\end{equation}
The crucial fact about singular cochains is that
the rows of the double complex \eqref{SingularDoubleComplex}
are acyclic except at the first column, where the horizontal cohomology
is equal to $C^{\bullet}(M, U)$. More precisely, 
consider the extended $p^{\mathrm{th}}$ row:
\begin{equation}\label{CpComplex}
0 \rightarrow
C^p ( M, U ) \xrightarrow{\upsilon^*}
C^p ( U ) \xrightarrow{\delta}
C^p ( U \times_M U ) \xrightarrow{\delta}
C^p ( U \times_M U \times_M U ) \xrightarrow{\delta}
\cdots
\end{equation}
Here the second arrow is the pull-back under the covering map 
$U \xrightarrow{\upsilon} M$. It follows from the
definition of $U \times_M U$ that 
\begin{equation}\label{DeltaUpsilon}
\delta \circ \upsilon^* = 0
\end{equation}
and hence \eqref{CpComplex} is a complex.
We claim that this complex is acyclic, 
in fact homotopic to $0$-complex. It is easier to describe
the dual homotopy for the dual complex of small singular chains 
\begin{equation}\label{CpDualComplex}
0 \leftarrow
C_p ( M, U ) \xleftarrow{\upsilon_{\ast}}
C_p ( U ) \xleftarrow{\delta^*}
C_p ( U \times_M U ) \xleftarrow{\delta^*}
C_p ( U \times_M U \times_M U ) \xleftarrow{\delta^*}
\cdots
\end{equation}
The following construction was explained to us by Matthew Ando.
Choose a section $\tau$ of the map $\upsilon_{\ast}$
($\tau$ exists since we consider $U$-small simplexes) and consider
a map from the total space of the complex
\eqref{CpDualComplex} to itself given by
\begin{align*}
\rho (\sigma) &= \tau (\sigma) \in C_p ( U ) & &\text{if }
\sigma \in C_p ( M, U ) \\
\rho (\sigma) &= 
\bigl( \tau \circ \upsilon_{\ast} \circ \mathrm{pr}_{\ast} ( \sigma ) 
\bigr) \times \sigma 
\in C_p ( \underbrace{U \times_M \ldots \times_M U}_{q+1} ) 
& &\text{if }
\sigma \in C_p ( \underbrace{U \times_M \ldots \times_M U}_{q} ) 
\end{align*}
where $\mathrm{pr} : U \times_M \ldots \times_M U \rightarrow U$ is
the projection map (onto any of the $U$-factors). One easily checks 
that 
\begin{align*}
\upsilon_{\ast} \circ \rho &= 
\mathrm{id} & &\text{on } C_p ( M, U ) \\
\rho \circ \upsilon^*  + \delta^* \circ \rho &= 
\mathrm{id} & &\text{on } C_p ( U ) \\
\rho \circ \delta^*  + \delta^* \circ \rho &= 
\mathrm{id} & &\text{on } 
C_p ( U \times_M \ldots \times_M U ) 
\end{align*}
which means that the complex \eqref{CpDualComplex}
(and thus the dual complex \eqref{CpComplex}) is homotopic to $0$.

\emph{Step 3.} 
Recall that $\epsilon$ is the inclusion of the first column in 
the double complex \eqref{SingularDoubleComplex}. 
Because of \eqref{DeltaUpsilon} the map $\epsilon \circ \upsilon^* :
C^{\bullet}(M ,U) \rightarrow
C(U)^{\bullet}_{\mathrm{tot}}$ is a morphism of complexes.
Step 1 reduces the Theorem to the statement that
the corresponding functor 
\begin{equation}\nonumber
\mathcal{H}^1 ( \epsilon \circ \upsilon^* ) : 
\mathcal{H}^1 ( C^{\bullet}(M ,U) ) \rightarrow
\mathcal{H}^1 ( C(U)^{\bullet}_{\mathrm{tot}} )
\end{equation}
is an equivalence of categories. Explicitly,
we have to show (a) the functor is full and faithful
and (b) any object in the range is isomorphic to an object 
in the image. Both statements follow from the acyclicity of 
the complex \eqref{CpComplex}:

(a) Let $x_1$, $x_2$ be objects in the image of 
$\mathcal{H}^1 ( \epsilon \circ \upsilon^* )$. 
This just means that $x_i = \upsilon^* (y_i) \in C^1 (U)$, 
where $y_i \in C^1 (M, U)$, $d y_i =0$.
Now we have 
\begin{multline}\nonumber
\mathrm{Hom}_{\mathcal{H}^1 ( C(U)^{\bullet}_{\mathrm{tot}} )} (x_1, x_2) = 
\{ b \in C^0 (U) \ | \  db = x_2 - x_1, \ \delta b = 0 \} = \\
= \{ b = \upsilon^* (a) \ | \ a \in C^0 (M, U)\ , \  da = y_2 - y_1 \} =
\mathrm{Hom}_{\mathcal{H}^1 ( C^{\bullet}(M,U) )} (y_1, y_2)
\end{multline}
where the second equality is due to acyclicity of 
the complex \eqref{CpComplex}
with $p=0$, $1$. 

(b) Consider an object of
$\mathcal{H}^1 ( C(U)^{\bullet}_{\mathrm{tot}} )$, i.e.
a pair \mbox{$(z,t) \in C^1 (U) \oplus C^0 ( U \times_M U )$,}
such that $\delta t = 0$, $dz = 0$, $dt = \delta z$.
Because the complex \eqref{CpComplex} is acyclic with $p=0$, there is
an element $c \in C^0(U)$ such that $\delta c = t$.
Now we think of $-c$ as an isomorphism in 
$\mathcal{H}^1 ( C(U)^{\bullet}_{\mathrm{tot}} )$ from
$(z,t)$ to $(z-dc, 0)$. Moreover, $\delta (z-dc) = 0$
and acyclicity of the complex \eqref{CpComplex} with $p=1$ implies
$(z-dc) = \epsilon \circ \upsilon^* (y)$ for some 
$y \in C^1 (M,U)$, $dy=0$. In other words,
$(z-dc, 0)$ is in the image of the functor 
$\mathcal{H}^1 ( \epsilon \circ \upsilon^* )$.

This completes the proof of the Proposition.
\end{proof}

\subsection{De Rham stack.}\label{DeRhamChains}
Similarly to singular cochains, one can consider the presheaf 
$\Omega^{\bullet}$ of differential forms with smooth coefficients.
\begin{proposition}One has:
\begin{enumerate}
\item\label{OmegaSheaf}
$\Omega^{p}$ is a sheaf for any $p$.
\item
$\mathcal{H}^0 ( \Omega^{\bullet} )$ is a sheaf (discrete stack).
\item
$\mathcal{H}^1 ( \Omega^{\bullet} )$ is a stack.
\end{enumerate}
\end{proposition}
\begin{proof}
The proof is similar to the proof of Proposition \ref{SingularChains}.
An important difference is that every differential form can
be pulled back along a covering map 
$U=\bigsqcup_{\alpha} U_{\alpha} \xrightarrow{\upsilon} M$ . 
Hence we don't need the first step of that proof 
in the current situation and, moreover,
we get the additional result \eqref{OmegaSheaf}. 
The rest of the proof is the same once one shows that Mayer-Vietoris
principle holds for $\Omega^p$, i.e. that the complex
\begin{equation}\label{OmegapComplex}
0 \rightarrow
\Omega^p ( M ) \xrightarrow{\upsilon^*}
\Omega^p ( U ) \xrightarrow{\delta}
\Omega^p ( U \times_M U ) \xrightarrow{\delta}
\Omega^p ( U \times_M U \times_M U ) \xrightarrow{\delta}
\cdots
\end{equation}
is acyclic. The proof of Mayer-Vietoris principle for
singular chains relied on the existence of a section of 
the push-forward map $\upsilon_*$. In the context of differential
forms the role of such a section is played by a partition of unity.
We recall the standard argument below.

First note that 
$\underbrace{U \times_M \ldots \times_M U}_{q} = 
\bigsqcup_{\alpha_1 \ldots \alpha_q} U_{\alpha_1} \cap \ldots \cap U_{\alpha_q}$. 
Given a form
\mbox{$\omega \in \Omega^p( \underbrace{U \times_M \ldots \times_M U}_{q})$}
we denote by $\omega_{\alpha_1 \ldots \alpha_q}$
its restriction to 
$U_{\alpha_1} \cap \ldots \cap U_{\alpha_q}$.
Now let $1= \sum_{\alpha} \eta_{\alpha}$ be a partition of unity
subordinate to the covering
$U=\bigsqcup_{\alpha} U_{\alpha} \xrightarrow{\upsilon} M$
and consider a map from the total space of the complex
\eqref{OmegapComplex} to itself given by
\begin{equation}\nonumber
\bigl( \rho ( \omega ) \bigr)_{\alpha_1 \ldots \alpha_q} = 
\sum_{ \alpha_0 } \eta_{\alpha_0} \omega_{\alpha_0 \alpha_1 \ldots \alpha_q}
\end{equation}
Then
\begin{align*}
\rho \circ \upsilon^{\ast} &= 
\mathrm{id} & &\text{on } \Omega^p ( M ) \\
\upsilon^* \circ \rho + \rho \circ \delta  &= 
\mathrm{id} & &\text{on } \Omega^p ( U ) \\
\delta \circ \rho + \rho \circ \delta &= 
\mathrm{id} & &\text{on } 
\Omega^p ( U \times_M \ldots \times_M U ) 
\end{align*}
which means that the complex \eqref{OmegapComplex}
is homotopic to $0$.
\end{proof}

\subsection{(Pre)sheaf (hyper)cohomology on stacks.}
Let $F^{\bullet}$ be a complex of pre-sheaves
and $\mathcal{W}$ be a stack. There are several ways to define
cohomology of $F^{\bullet}$ on $\mathcal{W}$. The most abstract
definition uses injective resolution, which is not very useful for
actual calculations. A more explicit construction (generalizing
\v{C}ech cohomology) exists for differentiable stacks. One starts with 
an atlas $\Gamma_{\bullet} \rightarrow \mathcal{W}$ and defines 
$H^n ( \mathcal{W} , F^{\bullet} )$ as the total cohomology of
the double complex $F^{\bullet} ( \Gamma_{\bullet} )$. It can be shown 
(see e.g.  \cite{Behrend2004}) that if
$F^{\bullet}$ is \v{C}ech-acyclic for covers of manifolds, then
$H^n ( \mathcal{W} , F^{\bullet} )$ does not depend on the 
choice of the atlas (for a general $F^{\bullet}$, 
the double complex will be the first term of a spectral sequence 
converging to $H^n ( \mathcal{W} , F^{\bullet} )$.

An intristic definition of hyper-cohomology uses the notion of
an Eilenberg-MacLane stack $\mathcal{K} ( F^{\bullet}, n )$, 
which is characterized by its universal
property (cf. \cite{Toen2006}): 
for any manifold $M$ the isomorphism classes of objects
of $(\mathcal{K} ( F^{\bullet}, n )) (M)$ are in (functorial in $M$)
bijection with $H^n ( F^{\bullet} (M) )$. 
Theorems \ref{SingularChains} and \ref{DeRhamChains} say that,
if $F^{\bullet}=C^{\bullet}$ or $\Omega^{\bullet}$, and $n = 1$ (or $0$), 
then $\mathcal{K} ( F^{\bullet}, n )$ is equivalent to 
$\mathcal{H}^n ( F^{\bullet} (M) )$ (or, even better, 
$\mathcal{H}^n ( F^{\bullet} (M) )$ provides an explicit 
construction of $\mathcal{K} ( F^{\bullet}, n )$)
Given an Eilenberg-MacLane stack $\mathcal{K} ( F^{\bullet}, n )$
the category $\mathcal{H}^n ( \mathcal{W} , F^{\bullet} ):=
\mathcal{H}om ( \mathcal{W}, \mathcal{K} ( F^{\bullet}, n ) )$
is called the cohomology category of the complex
$F^{\bullet}$ on the stack $\mathcal{W}$, and the set of
isomorphism classes of objects of this category the cohomology
$H^n ( \mathcal{W} , F^{\bullet} )$ of $F^{\bullet}$ on $\mathcal{W}$.
This definition assumes nothing about the stack $\mathcal{W}$.

If $n > 1$, then the presheaves of groupoids $\mathcal{H}^n (C^\bullet)$ 
and $\mathcal{H}^n (\Omega^\bullet)$ are not stacks
(and neither $\mathcal{K} ( C^{\bullet}, n )$ nor
$\mathcal{K} ( \Omega^{\bullet}, n )$ exists as a stack).
For example, the automorphism group of the $0$-object of 
$\mathcal{H}^2 (C^\bullet)$ is the first singular
cohomology group and so is not local. The reason is clear from the above 
proofs of locality: to make double complexes work
we need to consider cochains of all degrees lower than $n$ even if
we are only interested in $n^{\mathrm{th}}$ cohomology.
This leads one to think of, say, 
$H^n (C^{\bullet} (M))$ as the set of isomorphism classes of objects of 
an $n$-category having cochains of degree $k \leq n$ as $(n-k)$-morphisms.
So for each manifold we have a $n$-category \
$\mathcal{H}^n (C^{\bullet} (M))$, albeit 
a very simple one -- with invertible morphisms and strict associativity 
in all degrees, and this sheaf of $n$-categories satisfies descent 
condition. We refer the reader to 
\citelist{\cite{Toen2006} \cite{Lurie2006}} for the theory of
higher stacks. In particular, 
Theorems \ref{SingularChains} and \ref{DeRhamChains}
generalize to higher degrees: 
$\mathcal{K} ( C^{\bullet}, n )$ is equivalent to 
$\mathcal{H}^n ( C^{\bullet} (M) )$ and 
$\mathcal{K} ( \Omega^{\bullet}, n )$ is equivalent to 
$\mathcal{H}^n ( \Omega^{\bullet} (M) )$ for any $n$.

In the present paper we avoid higher stacks and hence higher cohomology.
There is a situation however in which one can describe the second
cohomology group of a complex of presheaves in terms of the usual 
(1-)stacks. This happens if $H^0 (M)$ is trivial for any $M$. In that case
the 2-categories involved are 2-discreet and thus equivalent to 
1-categories. An example of such situation is provided 
by differential characters.

\subsection{Differential characters.}\label{DCChains}

Recall that $C^\bullet$ is the smooth singular cochain complex. Let
$C_{\mathbb{R}}^\bullet = C^\bullet \otimes_{\mathbb{Z}} \mathbb{R}$ and,
for consistency, denote $C^\bullet$ by $C_{\mathbb{Z}}^\bullet$. 
The integration provides an inclusion of the de Rham complex
$\Omega^{\bullet}$ into  $C_{\mathbb{R}}^\bullet$.

Given an integer $s \geq 0$ Hopkins-Singer \cite{HopkinsSinger2005}
define the following complex of presheaves of abelian groups 
over $\mathcal{M}an$:
\begin{equation}\nonumber
DC^n_s = \{ ( c, h, \omega ) \ | \ \omega = 0 \text{ if } n < s \} 
\subset C^n_{\mathbb{Z}} \times C^{n-1}_{\mathbb{R}} \times \Omega^n 
\end{equation}
with the differential 
\begin{equation}\nonumber
d ( c, h, \omega ) = ( dc, \omega - c- dh, d\omega ) 
\end{equation}
The importance of this complex is that the cohomology group
$H^k (DC^\bullet_k (M))$ is isomorphic to the group of
Cheeger-Simons differential characters \cite{CheegerSimons1985} 
of degree $k$ on a manifold $M$. Recall that such a differential
character is a pair $(\omega , \chi )$, where $\omega \in \Omega^k (M)$ 
is a differential k-form and 
$\chi : Z_{k-1} (M) \rightarrow \mathbb{R}/\mathbb{Z}$ is a character 
of the group of smooth singular $(k-1)$-cycles. 
This pair should satisfy the following condition: 
\begin{equation}\label{OmegaChiCond}
\chi ( \partial S )  \equiv \int_{S} \omega 
\quad \mathrm{mod} \ \mathbb{Z} 
\end{equation}
for any smooth singular $k$-chain $S$.
The isomorphism between $H^k (DC^\bullet_k (M))$ and the group
of differential characters is given by the map 
$[(c,h,\omega)] \mapsto (\omega, \chi)$, 
where $\chi$ is the restriction of $h \in C^{k-1}_{\mathbb{R}} (M)$
to $Z_{k-1} (M)$ modulo $\mathbb{Z}$. 
The condition \eqref{OmegaChiCond} ensures that
$c = \omega - dh$ is an integral cochain.

We are going to show that $H^2 (DC^\bullet_1 )$ and 
$H^2 (DC^\bullet_2 )$ classify $S^1$-bundles and 
$S^1$-bundles with connections respectively. So in the cases
we are interested in we always have $s > 0$.

Following our usual procedure we replace the second cohomology 
with the presheaf of groupoids 
$\mathcal{DC}^2_s := \mathcal{H}^2 ( DC^\bullet_s )$.
At first sight $\mathcal{DC}^2_s$ has no chance of being 
a stack since we discarded $0$-cochains. Fortunately, if 
$s > 0$, then $DC^\bullet_s$ has trivial $0^{\mathrm{th}}$ cohomology because
\begin{equation}\nonumber
d (c, 0, 0) = (dc, -c, 0) 
\quad \text{ for } \quad 
(c, 0, 0) \in DC^0_s (M),  \ s>0
\end{equation}
This vanishing allows us to prove the following crucial result.
\begin{proposition}
If $s > 0$ then the presheaf of groupoids $\mathcal{DC}_s^2$ is a stack.
\end{proposition}
We call $\mathcal{DC}_2^2$ the stack of differential characters
(of degree 2).
\begin{proof}
Let us consider two double complexes associated to the complex 
$DC_s^{\bullet}$ and a covering $U \rightarrow M$. The first is 
the usual one computing the descent data
\begin{equation}\label{DCDoubleComplex}
\xymatrix{
& {\vdots} & {\vdots} & {\vdots} & \\
0 \ar[r] &
{DC_s^2 ( M )} 
\ar[r]^-{\upsilon^*} \ar[u]^{d} &
{DC_s^2 ( U )} 
\ar[r]^-{\delta} \ar[u]^{d} &
{DC_s^2 ( U \times_M U )} 
\ar[r]^-{\delta} \ar[u]^{-d} &
{\cdots}
\\
0 \ar[r] &
{DC_s^1 ( M )} 
\ar[r]^-{\upsilon^*} \ar[u]^{d} &
{DC_s^1 ( U )} 
\ar[r]^-{\delta} \ar[u]^{d} &
{DC_s^1 ( U \times_M U )} 
\ar[r]^-{\delta} \ar[u]^{-d} &
{\cdots}
\\
0 \ar[r] &
{DC_s^0 ( M )} 
\ar[r]^-{\upsilon^*} \ar[u]^{d} &
{DC_s^0 ( U )} 
\ar[r]^-{\delta} \ar[u]^{d} &
{DC_s^0 ( U \times_M U )} 
\ar[r]^-{\delta} \ar[u]^{-d} &
{\cdots}}
\end{equation}
and the second is
\begin{equation}\label{H1DCDoubleComplex}
\xymatrix{
& {\vdots} & {\vdots} & {\vdots} & \\
0 \ar[r] &
{DC_s^2 ( M )} 
\ar[r]^-{\upsilon^*} \ar[u]^{d} &
{DC_s^2 ( U )} 
\ar[r]^-{\delta} \ar[u]^{d} &
{DC_s^2 ( U \times_M U )} 
\ar[r]^-{\delta} \ar[u]^{-d} &
{\cdots}
\\
0 \ar[r] &
{\widetilde{DC_s}^1 ( M )} 
\ar[r]^-{\upsilon^*} \ar[u]^{d} &
{\widetilde{DC_s}^1 ( U )} 
\ar[r]^-{\delta} \ar[u]^{d} &
{\widetilde{DC_s}^1 ( U \times_M U )} 
\ar[r]^-{\delta} \ar[u]^{-d} &
{\cdots}}
\end{equation}
Here
$\widetilde{DC_s}^1 ( \cdot )= DC_s^1 ( \cdot )/ d \ DC_s^0 ( \cdot )$. 
We denote $[x] \in \widetilde{DC_s}^1 ( \cdot )$ the class of 
\mbox{$x \in DC_s^1 ( \cdot )$.}

Rows of \eqref{DCDoubleComplex} are acyclic since they are
direct sums of rows of the double complexes associated
to $C_{\mathbb{Z}}^{\bullet}$, $C_{\mathbb{R}}^{\bullet}$, and 
$\Omega^{\bullet}$. However in order to prove the Proposition we
need the acyclicity of the rows of \eqref{H1DCDoubleComplex}.
Only the last row is new. So let 
\mbox{$[x] \in \widetilde{DC_s}^1 ( \cdot )$}, $\delta [x] = [0]$.
This means that $\delta x = dt$ for some $t$. By anticommutativity 
of the double complex we have $d \delta t = - \delta^2 x = 0$. 
But columns of
\eqref{DCDoubleComplex} have trivial $0^{\mathrm{th}}$ cohomology.
Hence $\delta t =0$ and so, $0^{\mathrm{th}}$ row of  
\eqref{DCDoubleComplex} being acyclic, $t=\delta s$ for some $s$. 
It follows that
$\delta (x+ds) =0$ and the acyclicity of the $1^{\mathrm{st}}$ row
of \eqref{DCDoubleComplex} implies $x+ds = \delta y$ for some $y$.
Then $[x] = \delta [y]$, which shows that the last row of
\eqref{H1DCDoubleComplex} is acyclic. Of course one can replace
this diagram chasing with the long exact sequence in cohomology
associated to the short exact sequence 
$0 \rightarrow DC_s^0 ( \bullet ) \xrightarrow{d} DC_s^1 ( \bullet ) 
\rightarrow \widetilde{DC_s}^1 ( \bullet ) \rightarrow 0$ of rows
(i.e. $\delta$-complexes). However we prefer explicit calculations
with chains because they are easily translated into the language
of chain categories.

The rest of the proof is the same as for 
the singular and the de Rham complexes. 
\end{proof}

\subsection{Equivariant differential characters.}

Let $\Gamma_{\bullet}$ be a Lie groupoid. By definition (or rather, by
an argument similar to the proof of proposition \ref{H01Total}), 
the category 
$\bigl( \mathcal{H}^2 (DC_s^{\bullet} ) \bigr) (\Gamma_{\bullet} )$ 
is isomorphic
to the second cohomology category of the total complex of
\eqref{H1DCDoubleComplex}. However, for computational 
as well as aesthetic purposes, we would like to use the 
standard double complex $DC_s^{\bullet} (\Gamma_{\bullet})$. Fortunately
one has the following result.

\begin{proposition}If $s>0$ then
the following categories are equivalent
\begin{equation}
\bigl( \mathcal{H}^2 ( DC_s^\bullet ) \bigr) ( \Gamma_{\bullet} ) = 
\mathcal{H}^2 ( DC_s(\Gamma)_{\mathrm{tot}}^{\bullet} ) 
\end{equation}
\end{proposition}
\begin{proof}
The proof is a combination of the proofs of
Propositions \ref{H01Total} and \ref{DCChains}
and relies on the fact that $H^0 (DC_s^{\bullet} (M)) = 0$ 
for any manifold $M$.
\end{proof}

\section{Chern functor and prequantization.}
\label{ChernFunctor}

In this section we prove classification theorems for principal 
$S^1$-bundles with or without connections.

\subsection{Chern functor.}\label{ChernClassTheorem}
There is a natural functor 
$\mathrm{Ch}_{\mathrm{triv}} : \mathcal{BS}_{\mathrm{triv}}^1 \rightarrow 
\mathcal{DC}_1^2$ defined as follows.
\begin{itemize}
\item
The category $\mathcal{BS}_{\mathrm{triv}}^1 (M)$ has unique
object $M \times S^1$.
We put 
\begin{equation}\nonumber
\mathrm{Ch}_{\mathrm{triv}} ( M \times S^1 ) = 
\underline{0}=(0,0,0) \in DC_1^2 (M) \ .
\end{equation} 
\item
A morphism in $\mathcal{BS}_{\mathrm{triv}}^1 (M)$ is a smooth function
$f: M \rightarrow S^1 = \mathbb{R}/\mathbb{Z}$. 
We pick a lift $\tilde{f} : M \rightarrow \mathbb{R}$ ($\tilde{f}$
is not required to be smooth) and put
\begin{equation}\nonumber
\qquad \qquad
\mathrm{Ch}_{\mathrm{triv}} (f)= [(d(\tilde{f}-f), -\tilde{f}, -d f)] 
\in \mathrm{Hom}_{\mathcal{DC}_1^2 (M)} (\underline{0},\underline{0}) =
H^1 ( DC^{\bullet}_1 (M)) \ .
\end{equation}
The value $\mathrm{Ch}_{\mathrm{triv}} (f)$ does not depend 
on the choice of $\tilde{f}$.
\end{itemize}
If $M$ is contractible (for example $M=\mathbb{R}^n$) then
we have
\begin{equation}\nonumber
\begin{split}
H^2 ( DC^{\bullet}_1 (M)) &= 0 \ , \\
H^1 ( DC^{\bullet}_1 (M)) &= \Omega^0 (M)/ \{
\text{const functions } M \rightarrow \mathbb{Z} \} =
\{ \text{$C^{\infty}$ functions } M \rightarrow S^1 \} \ ,
\end{split}
\end{equation}
which means that $\mathrm{Ch}_{\mathrm{triv}}$ is an equivalence of categories 
from $\mathcal{BS}^1_{\mathrm{triv}} (\mathbb{R}^n)$ 
to $\mathcal{DC}_1^2 (\mathbb{R}^n)$.
Therefore (by the universal property of stackification)
$\mathrm{Ch}_{\mathrm{triv}}$ induces an equivalence 
$\mathrm{Ch}$ from the stackification of
$\mathcal{BS}^1_{\mathrm{triv}}$ (i.e. $\mathcal{BS}^1$) 
to the stack $\mathcal{DC}_1^2$. 
Note that the functor $\mathrm{Ch}$ is defined only up to 
a natural transformation.

To summarize we have the following theorem.
\begin{theorem}
There is an equivalence of stacks $\mathrm{Ch}: \mathcal{BS}^1
\rightarrow \mathcal{DC}_1^2$. In particular:
\begin{itemize}
\item
For any stack $\mathcal{W}$ the functor $\mathrm{Ch}$ induces 
an equivalence of categories from 
$\mathcal{BS}^1 (\mathcal{W})$ to $\mathcal{DC}_1^2 (\mathcal{W}) = 
(\mathcal{H}^2 ( DC_1^{\bullet})) (\mathcal{W})$;
\item
For any manifold $M$ the functor $\mathrm{Ch}$ induces a bijection from
the set of isomorphism classes of principal $S^1$-bundles on $M$ to
$H^2 ( DC_1^{\bullet} (M))$;
\item 
For any groupoid $\Gamma_{\bullet}$ the functor $\mathrm{Ch}$ 
induces a bijection from the set of isomorphism classes of 
$\Gamma_{\bullet}$-equivariant principal $S^1$-bundles to
the second total cohomology group of the complex
$DC_1^{\bullet} ( \Gamma_{\bullet})$;
\item 
For a Lie group $G$ acting on a manifold $M$ 
the functor $\mathrm{Ch}$ 
induces a bijection from the set of isomorphism classes of 
$G$-equivariant principal $S^1$-bundles on $M$ to
the equivariant cohomology group
$H^2_G (DC_1^{\bullet} ( M ))$.
\end{itemize}
\end{theorem}

\noindent
Here by equivariant cohomology we
mean the simplicial model of equivariant cohomology, that is
$H^2_G (DC_1^{\bullet} ( M )) := 
H^2 (DC_1^{\bullet} ( M \leftleftarrows G \times M ))$.
The reason for separating this case is that there are other
models (especially, for a compact group action), but of course
before using them one has to prove they give the same answer.

\subsection{Weil Theorem.} 
Let us explain how the above result implies Weil's classification theorem
for principal $S^1$-bundles on a manifold $M$.
Recall that 
\begin{equation}\nonumber
DC^n_1 (M) = \{ ( c, h, \omega ) \ | \ \omega = 0 \text{ if } n = 0 \} 
\subset C^n_{\mathbb{Z}} (M) \times C^{n-1}_{\mathbb{R}} (M) \times 
\Omega^n (M) \ . 
\end{equation}
We claim that the projection 
$p : DC^2_1 (M) \rightarrow C^2_{\mathbb{Z}} (M)$ 
induces an isomorphism of the second cohomology groups.
Both surjectivity and injectivity follow from the de Rham Theorem which 
says that the inclusion (given by integration) 
of the complex of differential forms into 
the complex of singular cochains with real coefficients
induces isomorphism in cohomology. In plain words it means
that, given $a \in C^n_{\mathbb{R}} (M)$ such that  $da \in \Omega^{n+1} (M)$, 
there exist $\alpha \in \Omega^n (M)$ and $b \in C^{n-1}_{\mathbb{R}} (M)$
such that $a = \alpha + db$. We call $\alpha$ a de Rham representative
of (the class) of $a$.

Let $\omega \in \Omega^2 (M)$ be a de Rham representative of a cocycle 
$c \in C^2_{\mathbb{Z}} (M) \subset C^2_{\mathbb{R}} (M)$, i.e. 
$c = \omega - dh$, $h \in C^1_{\mathbb{R}}$.
Then the cochain $(c,h,w)$ is closed in $DC^2_1 (M)$ and provides an extension 
of $c$. This proves surjectivity of $p$. To prove injectivity consider 
$x \in DC^2_1 (M)$ such that $dx = 0$, $p(x) = db$, or
$x=(db, h, d(b+h))$.
Let $\alpha \in \Omega^1 (M)$ be a de Rham representative of $b+h$,
i.e. $b + h = \alpha - df$, $f \in C^0_{\mathbb{R}}$. 
Then $x = d(b, f, \alpha)$, which implies that $p$ is
injective on cohomology. 

Composing the projection $p$ with the bijection from the manifold
case of Theorem \ref{ChernClassTheorem} we obtain Weil's theorem.

\begin{theorem}\label{WeilTheorem}
There is a bijection (the first Chern class)
form the set of isomorphism classes of principal $S^1$-bundles 
on a manifold $M$ to $H^2 ( M, \mathbb{Z} )$.
\end{theorem}

A few remarks are in order. 

First, one can consider the natural group structures on  
the sets $H^2 ( M, \mathbb{Z} )$, $H^2 ( DC^{\bullet}_1 (M))$,
and on the sets of objects of the categories
$\mathcal{BS}^1 (M)$ and $\mathcal{DC}^{2}_1 (M))$.
Then by a straightforward refinement of the above argument 
(which requires considering
sheaves of groupoids with additive structure on objects, their descent and 
stackifications, etc.) one can show that the first Chern class 
is a group homomorphism.

Second, the reader could wonder why do we deal with
the complicated complex $DC^{\bullet}_1$ when the final answer 
involves only singular cochains. The reason 
is that we would like to get a \emph{local}
proof of the Weil's theorem. In other words we consider the
equivalence of stacks $\mathrm{Ch}$ to be more
fundamental than the global bijection between isomorphism classes
of objects. This local equivalence makes no sense for
the complex $C^{\bullet}_{\mathbb{Z}}$ of singular cochains since 
$\mathcal{H}^2 (C^{\bullet}_{\mathbb{Z}})$ is not a stack
(it is a 2-stack).

Finally, there is a sheaf-theoretic proof of the Weil's theorem.
One identifies $S^1$-bundles with the first \v{C}ech cohomology
of the sheaf $\underline{{S}^1}$ of smooth $S^1$-valued functions and
then uses the short exact sequence of sheaves 
$0 \rightarrow \underline{\mathbb{Z}} \rightarrow \underline{\mathbb{R}} 
\rightarrow \underline{{S}^1} \rightarrow 0$ to show that 
$\check{H}^1 ( M, \underline{{S}^1} ) = 
\check{H}^2 ( M, \underline{\mathbb{Z}} ) = H^2 ( M, \mathbb{Z} )$.
Let us explain the relation of this argument with our approach through
cochain stacks. Note that
even though $\underline{{S}^1}$ is a sheaf (hence a stack), 
the way it is used is not local - one has to consider a covering of
the manifold to evaluate the first cohomology of $\underline{{S}^1}$. 
In order to replace the first cohomology group by the group of global 
sections one has to shift degree by one. So we think of 
$\underline{{S}^1}$ as morphisms in a category
rather than objects. The resultant category is 
$\mathcal{BS}^1_{\mathrm{triv}}$. Now we have to stackify it. The
answer is $\mathcal{DC}^2_1 \cong \mathcal{BS}^1$. The exact sequence 
$0 \rightarrow \underline{\mathbb{Z}} \rightarrow \underline{\mathbb{R}} 
\rightarrow \underline{{S}^1} \rightarrow 0$ is hidden in the proof
of the fact that $\mathcal{DC}_1^2$ is a stack.

\subsection{Proper stacks and equivariant Weil Theorem.}
Consider the (action) groupoid $M \leftleftarrows G \times M$
associated to an action of a Lie group $G$ on a manifold $M$.
If the group $G$ is compact the standard averaging argument shows that 
the complex $\Omega^0 (\Gamma_{\bullet}) = C^{\infty} (\Gamma_{\bullet})$ 
is acyclic except in degree 0 (the cohomology of this complex is called
differentiable cohomology of $\Gamma_{\bullet}$). 
The vanishing of higher differential cohomology 
remains true for an arbitrary proper Lie
groupoid (see \cite{Crainic2003}) with essentially the same proof. 
Recall that a groupoid is proper if the map  
$(s,t) : \Gamma_1 \rightarrow \Gamma_0 \times \Gamma_0$ is 
proper. A stack $\mathcal{W}$ is called proper if it has a proper atlas. 

\begin{theorem}
Let $\Gamma_{\bullet}$ be a proper Lie groupoid (or more generally,
a Lie groupoid with vanishing higher differentiable cohomology). 
There is a bijection (the equivariant first Chern class)
from the set of isomorphism classes of 
$\Gamma_{\bullet}$-equivariant principal $S^1$-bundles to
the equivariant cohomology group $H^2 ( \Gamma_{\bullet}, \mathbb{Z} )$.
\end{theorem}
\noindent Here $H^2 ( \Gamma_{\bullet}, \mathbb{Z} )$ is the second total
cohomology of the double complex 
$C_{\mathbb{Z}}^{\bullet} ( \Gamma_{\bullet})$ of smooth 
equivariant singular cochains.

In the case of a compact group action this theorem was proved by
Brylinski \cite{Brylinski2000} using a \v{C}ech-type argument.

\begin{proof}
The proof is similar to the proof of the Weil theorem in 
subsection \ref{WeilTheorem}. 
We have to show that the second total cohomology groups
of the double complexes $C_{\mathbb{Z}}^{\bullet} ( \Gamma_{\bullet})$ 
and $DC_1^{\bullet} ( \Gamma_{\bullet})$ are isomorphic. Let us write
down explicitly cochains of
$DC_1^{\bullet} ( \Gamma_{\bullet})$ of total degree $n \leq 2$:
\begin{equation}\label{DC1Chains}
\xymatrix{\\
\{ (c_1 , h_1 , \omega _1 ) \} \ar[r]^-{\delta} \ar[u]^{d} & \\ 
\{ (b_1 , g_1 , \alpha _1 ) \} \ar[r]^-{\delta} \ar[u]^{d} & 
\{ (c_2 , h_2 , \omega _2 ) \} \ar[r]^-{\delta} \ar[u]^{-d} & \\
\{ (a_1 , 0 , 0) \} \ar[r]^-{\delta} \ar[u]^{d} & 
\{ (b_2 , 0 , 0) \} \ar[r]^-{\delta} \ar[u]^{-d} & 
\{ (c_3 , 0 , 0) \} \ar[r]^-{\delta} \ar[u]^{d} &
}
\end{equation}
Here $a_i$, $b_i$, $c_i$, are singular cochains with integer coefficients,
$g_i$, $h_i$, are singular cochains with real coefficients, and
$\alpha_i$, $\omega_i$, are differential forms. The projection onto the
first element in each triple provides a map 
$p : (DC_1 ( \Gamma))_{\mathrm{tot}}^{\bullet} \rightarrow
(C_{\mathbb{Z}} ( \Gamma ))_{\mathrm{tot}}^{\bullet}$ and we want to show 
that $p$ induces an isomorphism on the second total cohomology. 
To prove surjectivity we have to extend a 
$d_{\mathrm{tot}}$-closed cochain \mbox{$(c_1, c_2, c_3)$} in 
$(C_{\mathbb{Z}} ( \Gamma ))_{\mathrm{tot}}^2$, to
a $d_{\mathrm{tot}}$-closed cochain in 
$(DC_1 ( \Gamma ))_{\mathrm{tot}}^2$. First we extend 
$c_1$ to $(c_1, h_1, \omega_1)$, $d (c_1, h_1, \omega_1) = 0$,
using the same argument as in subsection \ref{WeilTheorem}.
Then we have 
\begin{equation}\nonumber
d c_2 = \delta c_1 = \delta( \omega_1 -d h_1 ) = 
\delta \omega_1 - d \delta h_1 \ .
\end{equation}
Hence 
\begin{equation}\nonumber
d (c_2 + \delta h_1) = \delta \omega_1 \in \Omega^2 ( \Gamma_1 )
\end{equation}
and, by the de Rham Theorem, we can find 
$h'_2 \in C^0_{\mathbb{R}} (\Gamma_1)$ and
$\omega'_2 \in \Omega^1 ( \Gamma_1 )$ such that
\begin{equation}
c_2 + \delta h_1 = \omega'_2 - d h'_2 
\end{equation}
or, in other words, $d (c_2, h'_2, \omega'_2) = \delta (c_1, h_1, \omega_1)$.
A similar argument shows that $c_3 - \delta h'_2 \in \Omega^0 (\Gamma_2)$.
But then, by the Theorem assumption, $c_3 - \delta h'_2 = \delta f$ 
for some $f \in \Omega^0 (\Gamma_1)$. Finally we put $h_2 = h_2' + f$,
$\omega_2 = \omega'_2 + df$ to get a $d_{\mathrm{tot}}$-closed cochain 
$((c_1,h_1,\omega_1), (c_2,h_2,\omega_2), (c_3,0,0))$ in
$(DC_1 ( \Gamma ))_{\mathrm{tot}}^2$. This proves surjectivity of $p$.

To prove injectivity of $p$ (on cohomology) 
we have to show (cf. \ref{WeilTheorem})
that any cocycle $x$ in $(DC_1 ( \Gamma ))_{\mathrm{tot}}^2$ of the form
$x=((db_1 , h_1 , \omega_1), (\delta b_1 - db_2 , h_2 , \omega_2), 
(\delta b_2 ,0,0))$ is in the image of
$d_{\mathrm{tot}}$. Repeating the argument in subsection
\ref{WeilTheorem} we can assume that $b_1=h_1 = \omega_1 = 0$.
Then $d_{\mathrm{tot}} x=0$ implies 
$d (h_2 - b_2) = \omega_2$ and $\delta (h_2 - b_2) =0$. 
Hence $(h_2 - b_2 ) \in \Omega^0 ( \Gamma_1 )$
and $\delta (h_2 -b_2) = 0$. Therefore, by the Theorem assumption,
$h_2 - b_2 = \delta f$ for some $f \in \Omega^0 (\Gamma_0 )$. This means 
$((0,0,0), (- db_2 , h_2 , \omega_2), (\delta b_2 , 0 , 0)) = 
d_{\mathrm{tot}} ((0,f,df), (b_2,0,0))$, which proves the injectivity of
$p$ on cohomology and hence the Theorem.
\end{proof}

\subsection{Prequantization.}
We now turn to bundles with connections and consider a functor 
$\mathrm{DCh}_{\mathrm{triv}} : 
\mathcal{DBS}_{\mathrm{triv}}^1 \rightarrow \mathcal{DC}_2^2$
given by (cf. \ref{ChernClassTheorem}):
\begin{itemize}
\item
$\mathrm{DCh}_{\mathrm{triv}} ((M \times S^1 , a + d\theta )) =
(0, a, da)$ on objects;
\item
$\mathrm{DCh}_{\mathrm{triv}} (f)= [(d(\tilde{f}-f), -\tilde{f}, 0)] 
\in DC^1_2 / d (DC^0_2)$ on morphisms.
Here we think of a morphism in $\mathcal{DBS}_{\mathrm{triv}}^1$ from 
$(M \times S^1 , a + d\theta )$ to $(M \times S^1 , a' + d\theta )$
as a smooth function
$f: M \rightarrow S^1 = \mathbb{R}/\mathbb{Z}$ such that 
$df = a' - a$ and let $\tilde{f} \in C^0_{\mathbb{R}} (M)$ be 
a lift of $f$. 
\end{itemize}

As in \ref{ChernClassTheorem}, a simple cohomology calculation in
$DC^{\bullet}_2 ( \mathbb{R}^n )$ shows that $\mathrm{DCh}_{\mathrm{triv}}$
restricts to an equivalence of categories between
$\mathcal{DBS}_{\mathrm{triv}}^1 ( \mathbb{R}^n )$ and 
$\mathcal{DC}_2^2 ( \mathbb{R}^n )$, and hence induces 
an equivalence of stacks $\mathrm{DCh}: \mathcal{DBS}^1 
\rightarrow \mathcal{DC}_2^2$.
The quasi-inverse functor $\mathrm{Preq}$ 
(defined up to a natural transformation)
is called the prequantization. 

To summarize we have the following theorem.

\begin{QTheorem}\label{PreQTheorem}
There is an equivalence of stacks 
\mbox{$\mathrm{Preq}: \mathcal{DC}_2^2 \rightarrow \mathcal{DBS}^1$} 
from the stack
of differential characters to the stack of principal $S^1$-bundles with 
connections. In particular, isomorphism classes of principal $S^1$-bundles 
with connections are classified by differential characters.
More precisely:
\begin{itemize}
\item
For any stack $\mathcal{W}$ the functor $\mathrm{Preq}$ induces 
an equivalence of categories from 
$\mathcal{DC}_2^2 (\mathcal{W}) = 
(\mathcal{H}^2 ( DC_2^{\bullet})) (\mathcal{W})$ to
$\mathcal{DBS}^1 (\mathcal{W})$;
\item
For any manifold $M$ the functor $\mathrm{Preq}$ induces a bijection from
$H^2 ( DC_2^{\bullet} (M))$ to the set isomorphism classes of 
principal $S^1$-bundles with connections on $M$;
\item 
For any groupoid $\Gamma_{\bullet}$ the functor $\mathrm{Preq}$ 
induces a bijection from 
the second total cohomology group of the double complex
$DC_2^{\bullet} ( \Gamma_{\bullet})$
to the set of isomorphism classes of 
$\Gamma_{\bullet}$-equivariant principal $S^1$-bundles with basic connections;
\item 
For a Lie group $G$ acting on a manifold $M$ 
the functor $\mathrm{Ch}$ induces a bijection from 
the equivariant cohomology group
$H^2_G (DC_2^{\bullet} ( M ))$ to 
the set of isomorphism classes of 
$G$-equivariant principal $S^1$-bundles with $G$-basic 
connections on $M$.
\end{itemize}
\end{QTheorem}

\noindent
To make things a bit more explicit let us write down cochains of
$DC_2^{\bullet} ( \Gamma_{\bullet})$ of total degree $n \leq 2$
on a groupoid $\Gamma_{\bullet}$
(cf. \eqref{DC1Chains} for cochains of $DC_1^{\bullet} ( \Gamma_{\bullet})$):
\begin{equation}\nonumber
\xymatrix{\\
\{ (c_1 , h_1 , \omega _1 ) \} \ar[r]^-{\delta} \ar[u]^{d} & \\ 
\{ (b_1 , g_1 , 0 ) \} \ar[r]^-{\delta} \ar[u]^{d} & 
\{ (c_2 , h_2 , 0 ) \} \ar[r]^-{\delta} \ar[u]^{-d} & \\
\{ (a_1 , 0 , 0) \} \ar[r]^-{\delta} \ar[u]^{d} & 
\{ (b_2 , 0 , 0) \} \ar[r]^-{\delta} \ar[u]^{-d} & 
\{ (c_3 , 0 , 0) \} \ar[r]^-{\delta} \ar[u]^{d} &
}
\end{equation}
Here $a_i$, $b_i$, $c_i$, are singular cochains with integer coefficients,
$g_i$, $h_i$, are singular cochains with real coefficients, and
$\omega_1$, is a differential form.
As in the manifold case (cf. subsection \ref{DCChains} and the Introduction), 
the map 
\begin{equation}\nonumber
\bigl[ \bigl( (c_1 ,h_1 ,\omega_1 ), 
(c_2, h_2 , 0), (c_3 , 0, 0) \bigr) \bigr] \mapsto
\bigl( \omega_1 , 
(h_1 \ \mathrm{mod} \, \mathbb{Z}, \, h_2 \ \mathrm{mod} \, \mathbb{Z})
\bigr)
\end{equation}
provides an isomorphism between
the second cohomology group 
$H^2_{\mathrm{tot}} ( DC_2^{\bullet} (\Gamma_{\bullet}))$ 
and the group of differential characters, i.e.
pairs $(\omega , \chi )$, where 
$\omega \in \Omega^2 (\Gamma_0 ) \subset \Omega^2 ( \Gamma_{\bullet} )$
and $\chi : Z_1 (\Gamma_{\bullet})
\rightarrow \mathbb{R}/\mathbb{Z}$ is a character 
of the group of smooth singular 1-cycles on the groupoid
$\Gamma_{\bullet}$, satisfying a condition similar 
(and generalizing) \eqref{OmegaChiCond}.
We refer the reader to \cite{LermanMalkinReduction} (where
a more general case is considered) for details, pictures, 
and examples.

\subsection{Equivariant Kostant theorem.}
Let $\Gamma_{\bullet}$ be a Lie groupoid.
We denote by  
$\Omega^2_{\mathbb{Z}, \mathrm{cl}, \mathrm{bas}} (\Gamma_{\bullet})$
the group of closed 2-forms on $\Gamma_0$ which are integral
(i.e. have integral periods, or equivalently, represent the
image of an integral cohomology class in de Rham cohomology),
and basic (i.e. are in the kernel of 
$\partial: \Omega^2 (\Gamma_0 ) \rightarrow \Omega^2 ( \Gamma_1 ))$.
Suppose 
$\omega_1 \in \Omega^2_{\mathbb{Z}, \mathrm{cl}, \mathrm{bas}} (\Gamma_{\bullet})$.
Then $(\omega_1, 0, 0) \in 
\Omega^2 (\Gamma_0) \times \Omega^1 (\Gamma_1) \times \Omega^0 (\Gamma_2)$
represents the image of a second integral cohomology class on $\Gamma_{\bullet}$ 
in the de Rham cohomology (here we used de Rham theorem on 
$\Gamma_{\bullet} \, $, cf. \cite{Behrend2004}, to identify
the de Rham and the real-valued singular cohomology). 
In other words, there exists an integer-valued 2-cocycle
\mbox{$(c_1 , c_2 , c_3 ) \in  
C_{\mathbb{Z}}^2 (\Gamma_0) \times C_{\mathbb{Z}}^1 (\Gamma_1) 
\times C_{\mathbb{Z}}^0 (\Gamma_2)$} 
and a real-valued 1-cochain 
\mbox{$(h_1 , h_2 ) \in 
C_{\mathbb{R}}^1 (\Gamma_0) \times C_{\mathbb{R}}^0 (\Gamma_1)$}
such that 
$d_{\mathrm{tot}} ( h_1, h_2 )  = (\omega_1, 0, 0)  - 
(c_1 , c_2 , c_3 )$, where $d_{\mathrm{tot}}$ is the total differential
in the double complex of smooth singular chains on $\Gamma_{\bullet}$.
This means that the projection map
\begin{gather}\nonumber
\eta: 
H^2_{\mathrm{tot}} (DC_2^{\bullet} ( \Gamma_{\bullet})) \rightarrow 
\Omega^2_{\mathbb{Z}, \mathrm{cl}, \mathrm{bas}} (\Gamma_{\bullet})
\\\nonumber
\eta \ ( \, [ ( (c_1 ,h_1 ,\omega_1 ), 
(c_2, h_2 , 0), (c_3 , 0, 0) )] \, ) = \omega_1
\end{gather}
is surjective, and it is easy to see that the kernel of $\eta$ is
isomorphic 
(by the map $
\bigl[ \bigl( (c_1 ,h_1 , 0 ), 
(c_2, h_2 , 0), (c_3 , 0, 0) \bigr) \bigr] \mapsto
[(h_1 \ \mathrm{mod} \, \mathbb{Z}, \, h_2 \ \mathrm{mod} \, \mathbb{Z})]$
to the cohomology group 
$H^1_{\mathrm{tot}} (C^\bullet_{\mathbb{R}/\mathbb{Z}} (\Gamma_{\bullet}))$ 
of smooth singular 1-cochains with values in 
$\mathbb{R}/\mathbb{Z}$. Hence we obtain a short exact sequence
\begin{equation}\label{KostantSequence}
0
\rightarrow 
H^1_{\mathrm{tot}} (C^\bullet_{\mathbb{R}/\mathbb{Z}} (\Gamma_{\bullet}))
\rightarrow
H^2_{\mathrm{tot}} (DC_2^{\bullet} ( \Gamma_{\bullet}) )
\rightarrow 
\Omega^2_{\mathbb{Z}, \mathrm{cl}, \mathrm{bas}} (\Gamma_{\bullet})
\rightarrow
0 \ .
\end{equation}
In the manifold case this sequence appeared in the original
Cheeger and Simons paper \cite{CheegerSimons1985}.
Since $H^2_{\mathrm{tot}} (DC_2^{\bullet} ( \Gamma_{\bullet}))$
classifies $\Gamma_{\bullet}$-equivariant principal $S^1$-bundles
with basic connections (cf. Theorem \ref{PreQTheorem}), 
we can interpret \eqref{KostantSequence} as follows:

\begin{theorem}
Let $\Gamma_{\bullet}$ be a Lie groupoid.
Given an integral closed basic $2$-form $\omega$ on $\Gamma_0$ there exists
a $\Gamma_{\bullet}$-equivariant principal 
$S^1$-bundle with a basic connection 
whose curvature is $\omega$. Moreover, the set of isomorphism 
classes of such bundles with connections is in bijection with 
$H^1_{\mathrm{tot}} (C^\bullet_{\mathbb{R}/\mathbb{Z}} (\Gamma_{\bullet}))$.
\end{theorem}

This theorem is due to Kostant \cite{Kostant1970} in the manifold case.
A more general version (about equivariant bundles with arbitrary connections)
was proved in \cite{BehrendXu2006} for proper groupoids, and in
\cite{LermanMalkinReduction} for arbitrary groupoids. 

As explained in the Introduction and remarks after Theorem \ref{WeilTheorem},
we believe that the abstract version
of Theorem \ref{PreQTheorem} is more fundamental than Cheeger-Simons- and
Kostant-type interpretations because (1) it describes the category of
bundles with connections (rather than isomorphism classes of objects) and
(2) it is a local statement. On the other hand,  
a concrete description in terms of cohomology groups 
or differential characters is obviously useful in calculations and/or
explicit geometric interpretations.

\end{document}